\documentclass[12pt]{amsproc}

\usepackage{amsmath}
\usepackage{amsthm}
\usepackage{amscd}
\usepackage{amssymb}
\usepackage{url}
\usepackage{nicefrac}
\usepackage[pdftex]{graphicx}%
\usepackage[outercaption]{sidecap}
\usepackage{amsfonts,epigraph}%
\usepackage{booktabs} 
\usepackage[all,cmtip]{xy}
\newtheorem{theorem}{Theorem}[section]
\newtheorem{corollary}[theorem]{Corollary}
\newtheorem{lemma}[theorem]{Lemma}
\newtheorem{proposition}[theorem]{Proposition}

\theoremstyle{definition}
\newtheorem{definition}[theorem]{Definition}

\newtheorem{choices-notations}[theorem]{Choices and Notations}

\newtheorem{rem}[theorem]{Remark}

\theoremstyle{remark}

\newtheorem{note}{Note\!\!}

\newcommand{\exendproof}{\renewcommand{\qed}{\relax}\end{proof}}

\newsavebox{\SmallMathBox}

\DeclareRobustCommand*{\nicefrac}[2]{\ifmmode\mathnicefrac{#1}
{ #2}%
  \else\textnicefrac{#1}{#2}\fi}
\newcommand*{\textnicefrac}[2]{\check@mathfonts%
\mbox{\raisebox{.5ex}{\fontsize\sf@size\z@\selectfont#1}\kern-.
1em%
/\kern-.1em\raisebox{- .25ex}{\fontsize\sf@size\z@\selectfont#2} }}
\newcommand*{\mathnicefrac}[2]{%
  \mathchoice
    {\m@fr@c{\scriptstyle}{#1}{#2}}
    {\m@fr@c{\scriptstyle}{#1}{#2}}
    {\m@fr@c{\scriptscriptstyle}{#1}{#2}}
    {\m@fr@c{\scriptscriptstyle}{#1}{#2}}}

\def\noi{\noindent}

\def\sqm1{\sqrt{-1}}

\def\tand{\mbox{\ \rm  and }}

\def\wh{\widehat}

\def\={\cong}
\def\>{\supset}
\def\<{\subset}

\def\12{\frac{1}{2}}
\def\0{^{\circ}}

\def\CC{{\mathbb C}}

\def\KK{{\mathbb K}}

\def\NN{{\mathbb N}}

\def\RR{{\mathbb R}}

\def\ZZ{{\mathbb Z}}

\def\Bb{{\mathcal B}}

\def\Ff{{\mathcal F}}

\def\Kk{{\mathcal K}}

\def\Ss{{\mathcal S}}

\def\C{\CC}

\def\e{\varepsilon}

\def\la{\lambda}

\def\N{\NN}

\def\R{\RR}

\def\Z{\ZZ}

\DeclareMathOperator{\defi}{def}
 \DeclareMathOperator{\dist}{dist}

\DeclareMathOperator{\dom}{dom}

  \DeclareMathOperator{\Graph}{graph}
 
 \DeclareMathOperator{\image}{im}
\DeclareMathOperator{\Index}{index}

\DeclareMathOperator{\nuli}{nul} \DeclareMathOperator{\ran}{im} 
  
 \DeclareMathOperator{\romS}{S} 

  \DeclareMathOperator{\Span}{span}

\begin{document}

\title{Gaps and relative dimensions}
%
%
%
%


\author{Chenfeng Liao}
\address{School of Mathematics Science, Peking University,
	Beijing 100871, P. R. China} \email{chenfengliao642@foxmail.com}

\author{Chaofeng Zhu}
\address{Chern Institute of Mathematics and Laboratory of Pure Mathematics and Combinatorics (LPMC), Nankai University,
	Tianjin 300071, P. R. China, https://orcid.org/0000-0003-4600-4253} \email{zhucf@nankai.edu.cn}
\thanks{Corresponding author: C. Zhu [\texttt{zhucf@nankai.edu.cn}]\\	
	Chaofeng Zhu is supported by National Key R\&D Program of China (2020YFA0713300), NSFC Grants (12571066, 11971245), Nankai Zhide Foundation and Nankai University.}

\date{today}

\subjclass[2010]{Primary 53D12; Secondary 58J30}

\keywords{gap between closed linear subspaces, semi-compact perturbation, relative dimension}


\begin{abstract}
In this paper, the notion of semi-compact perturbation of a closed linear subspace is introduced. Then for a of pair of closed linear subspace of a Banach space such that one is a semi-compact perturbation of the other, it is proved that the relative dimension between them is well-defined. If the perturbation is global, the relative dimension is stable, even the perturbed pair is a semi-compact perturbed one. After that, the notion of Fredholm tuple of closed linear subspaces in a Banach space is introduced. Then the stability of the Fredholm tuple is proved. Finally the perturbed augmented Morse index is studied. 
\end{abstract}

\maketitle


\section{Introduction}\label{l:introduction-semi-compact}

The virtual codimension between two pseudo-differential projections with the same principal symbol was introduced by L. G. Brown,
R. G. Douglas and P. A. Fillmore \cite{BrDF1973}, is a nice tool of measuring the difference between the dimensions of the ranges of the two projections. It behaves nicely and have many applications (for example, see \cite[Proposition 1.3]{DaiZh98}, \cite[Definition 7.1.5]{Lo02}). A natural question is:

{\em Can we define define the relative dimension between two closed linear subspaces in a Banach space when one is a compact perturbation of another, and develop nice theory on it?}  

The main purpose of the paper is to give an affirmative answer to this question. 

Let $X$ be a Banach space with two closed linear subspaces $M$ and $N$. If there exists a compact linear operator $K\in\Bb(X) $ such that $(I+K)M\subset N$, we call $N$ a {\em right global semi-compact perturbation} of $M$ and denote it by $M\sim^{gsc}N$. We call $N$ a {\em right semi-compact perturbation} of $M$ and denote it by $M\sim^{sc}N$ if $K$ is only required to be defined on $M$. In this case, we call the Fredholm index $\Index(I+K)\colon M\to N$ the {\em relative dimension} $[M-N]\in\Z\cup\{-\infty\}$ between $M$ and $N$. We define $M\overset{gsc}\sim N$ and call $N$ a {\em global semi-compact perturbation} of $M$ if either $M\sim^{gsc}N$ or $N\sim^{gsc}M$. Here we denote by $\ZZ$ the set of integers. For the details, see Definition \ref{d:semi-compact-perturtation} below.

We have the following stability theorem of relative dimension.

\begin{theorem}\label{t:stable-relative-dimention}
	Let $X$ be a Banach space with two closed linear subspaces $M$ and $N$. Assume that $N$ is a right global semi-compact perturbation of $M$. Then the following hold.
	\begin{itemize}
		\item[(a)] The function $[M-N]$ is well-defined.
		\item[(b)] Let $m$ be a positive integer. Then there are a compact linear operator $K_1\in\Bb(X)$ and two finite dimensional linear subspaces $V$ of $(I+K_1)M$, $U$ of $N$ respectively such that 
	    \begin{align}\label{e:finite-change-same-space}
	    	(I+K_1)M\oplus U&\begin{cases}
	    		=N\oplus V &\text{ if }[M-N]\in\ZZ,\\
	    		\subset N\oplus V&\text{ if }[M-N]=-\infty,
	    	\end{cases}\\
    	    \label{e:finite-change-same-space-dimension}
    	    \dim U&=\dim U+m \text{ if }[M-N]=-\infty,\\    	     \label{e:invertible-I+K1}
    	    (I+K_1)^{-1}&\in\Bb(X).
	    \end{align} 
        We set
        \begin{align*}
        	a\ :&=\|I+K_1\|\|(I+K_1)^{-1}\|.
        \end{align*}       
		\item[(c)] Let $m$ be a positive integer.
		Let $(M^\prime, N^\prime)$ be a pair of closed linear subspaces such that $N^\prime$ is a semi-compact perturbation of $M^\prime$.
		Set 
		\begin{align*}
			\theta_1\ :&=\frac{a(1+\gamma(U,(I+K_1)M))\delta(M^\prime,M)}{\gamma(U,(I+K_1)M)-a\delta(M^\prime,M)},\\
			\eta_1\ :&=(1+\delta(N,N^\prime)\gamma(N,V)^{-1})(1+\theta_1).
		\end{align*} 
	    Assume that 
	    \begin{align}\label{e:condition-relative-dimension-small}
	    	\delta(M^\prime, M)<a^{-1}\gamma(U,(I+K_1)M), \quad \eta_1<2.
	    \end{align}
	    Then we have  
		\begin{align}\label{e:stable-relative-dimention-small}
			[M^\prime-N^\prime]\le
			\begin{cases}
				[M-N] &\text{if }[M-N]\in\ZZ,\\
				-m &\text{if }[M-N]=-\infty.
			\end{cases}
		\end{align}
	    \item[(d)] Let $(M^\prime, N^\prime)$ be a pair of closed linear subspaces $N^\prime$ is a semi-compact perturbation of $M^\prime$ .
	    Set 
	    \begin{align*}
	    	\theta_2\ :&=\frac{a(1+\gamma(U,(I+K_1)N))\delta(N^\prime,N)}{\gamma(U,(I+K_1)N)-a\delta(N^\prime,N)},\\
	    	\eta_2\ :&=(1+\delta(M,M^\prime)\gamma(M,V)^{-1})(1+\theta_2).
	    \end{align*} 
	    Assume that 
	    \begin{align}\label{e:condition-relative-dimension-large}
	    	\delta(N^\prime, N)<a^{-1}\gamma(U,(I+K_1)M), \quad \eta_2<2.
	    \end{align}
	    Then we have   
	    \begin{align}\label{e:stable-relative-dimention-large}
	    	[M^\prime-N^\prime]\ge[M-N].
	    \end{align}
        \item[(e)] Assume that $[M-N]$ is an integer, i.e. $M\sim^{gc} N$. Let $(M^\prime, N^\prime)$ be a pair of closed linear subspaces with $M^\prime\overset{sc}\sim N^\prime$.
        Assume that both \eqref{e:condition-relative-dimension-small} and \eqref{e:condition-relative-dimension-large} hold.
        Then we have   
        \begin{align}\label{e:stable-relative-dimention}
        	[M^\prime-N^\prime]=[M-N].
        \end{align}		    
	\end{itemize} 
\end{theorem}

Even when $\dom(K^\prime)=X$, the lack of control of the norm $\|K^\prime-K\|$, makes the proof of Theorem \ref{t:stable-relative-dimention} rather complicated, where $K$ and $K^\prime$ are compact operators such that $(I+K)M\subset N$ and $(I+K^\prime)M^\prime\subset N^\prime$. Note that by \cite[Theorem 5(iv)]{Mu07}, for each compact linear operator $K$, for each $\varepsilon>0$, there exists a closed linear subspace $C$ with finite codimension such that $\|K|_C\|<\varepsilon$. This fact is not enough for us to prove Theorem \ref{t:stable-relative-dimention} in the special case when $\dom(K^\prime)=X$. 

The complemented case (Theorem \ref{t:stable-semi-compact-perturbations-complemented} below) give us evidence the general case. To overcome this difficulty, we firstly embed $X$ into the continuous function space $C(\Kk)$ by a linear isometric $j$, where $\Kk$ is a compact Hausdorff space (see Lemma \ref{l:embed-into-AP} below). A. Grothendieck \cite[Chapter I, p. 185]{Groth55} showed that $C(\Kk)$ has (MAP) (see definition \ref{d:Banach-AP} below). Then for a compact operator $T\in\Bb(\dom(T),X)$ with closed $\dom(T)$, the composite $jT\in\Bb(\dom(T),X)$ can be approximated by finite rank operators. This observation allows us proof Theorem \ref{t:semi-compact-perturbation} in whole generality.

We have the following generalization of \cite[Theorem]{HeKa70}.

\begin{theorem}\label{t:stable-relative-dimention-family}
	Let $X$ be a Banach space with a continuous family of pairs $\{(M(s),N(s));\}_{s\in J}$ of closed linear subspaces such that $N(s)$ is a global semi-compact perturbation of $M(s)$ for each $s\in J$, where $J$ is a connected topological space. Then the function $s\mapsto[M(s)-N(s)]$ is a constant function on $s\in J$. 	
\end{theorem}

Let $X$ be a Banach space with a pair of closed linear subspaces $(M,N)$. Let $\delta(M,N)$ and $\hat\delta(M,N)$ the gaps defined by \ref{d:closed-distance}.b below. By \cite[Corollary IV.2.6]{Ka95}, $\delta(M,N)<1$ implies $\dim M<\dim N$, and $\hat\delta(M,N)<1$ implies $\dim M=\dim N$. Then we may ask what will happened if we replace the inequalities conditions on the dimensions by the sufficiently small gaps. The classical stability theorem for Fredholm indices (\cite[Theorem IV.4.30]{Ka95}) is one of the results in this direction.

We have the following stability theorems for the indices of Fredholm tuples (cf. \cite[Theorem IV.4.24 and Theorem IV.4.30]{Ka95}). For the notion of index of semi-Fredholm tuple,
see Definition \ref{d:index-semi-Fredholm-tuple} below.

\begin{theorem}\label{t:stable-semi-Fredholm}
	Let $X$ be a Banach space with closed linear subspaces $Y_1$, $M$, $N$, $Y_2$ such that $M+N$ is closed and $Y_1\subset M\cap N\subset M+N\subset Y_2$ holds. Then the following hold.
	\begin{itemize}
		\item[(a)] Assume that $\dim Y_2/(M+N)$ is finite. Then for each tuple $(M^\prime, N^\prime;Y_2^\prime)$ of closed linear subspaces with $M^\prime+ N^\prime\subset Y_2^\prime$, and sufficiently small $\delta(M,M^\prime)$,
		$\delta(N,N^\prime)$ and $\delta(Y_2^\prime,Y_2)$, the space $M^\prime+N^\prime$ is closed, and we have
		\begin{align}\label{e:stable-sum}
			\dim Y_2^\prime/(M^\prime+ N^\prime)\le\dim Y_2/(M+N).
		\end{align}	
		\item[(b)] Assume that $\dim (M\cap N)/Y_1$ is finite. Then for each tuple $(Y_1^\prime;M^\prime, N^\prime)$ of closed linear subspaces with $Y_1^\prime\subset M^\prime\cap N^\prime$, and sufficiently small $\delta(Y_1,Y_1^\prime)$, $\delta(M^\prime,M)$ and 
		$\delta(N^\prime,N)$, the space $M^\prime+N^\prime$ is closed, and we have  
		\begin{align}\label{e:stable-cap}
			\dim (M^\prime\cap N^\prime)/Y_1^\prime\le\dim (M\cap N)/Y_1.		
		\end{align}			    
	\end{itemize}    
\end{theorem}

The most difficult step in the proof of Theorem \ref{t:stable-semi-Fredholm} is to prove the closeness of $M^\prime+N^\prime$. We will prove a splitting theorem (see Theorem \ref{t:existence-finite-dimension-embeding-close-oposite} below) to overcome the difficulty. Then we prove the closeness of $M^\prime+N^\prime$ by using the splitting theorem.  

\begin{theorem}\label{t:stable-Fredholm-index}
	Let $X$ be a Banach space with a semi-Fredholm tuple $(Y_1;M,N;Y_2)$ of closed linear subspaces. Then the following hold.
	\begin{itemize}		
		\item[(a)] Assume that the tuple  $(Y_1;M,N;Y_2)$ is Fredholm. Then for each tuple $(Y_1^\prime;M^\prime, N^\prime;Y_2^\prime)$ of closed linear subspaces with  $Y_1^\prime\subset M^\prime\cap N^\prime\subset M^\prime+N^\prime\subset Y_2^\prime$, and sufficiently small $\delta(Y_1^\prime,Y_1)$, $\delta(M,M^\prime)$,
		$\delta(N,N^\prime)$ and $\delta(Y_2^\prime,Y_2)$, we have  
		\begin{align}\label{e:stable-Fredholm-index-2}
			\Index(Y_1^\prime;M^\prime, N^\prime;Y_2^\prime)
			\ge	\Index(Y_1;M,N;Y_2).		
		\end{align}	
		\item[(b)] Assume that the tuple  $(Y_1;M,N;Y_2)$ is Fredholm. Then for each tuple $(Y_1^\prime;M^\prime, N^\prime;Y_2^\prime)$ of closed linear subspaces with  $Y_1^\prime\subset M^\prime\cap N^\prime\subset M^\prime+N^\prime\subset Y_2^\prime$, and sufficiently small $\delta(Y_1,Y_1^\prime)$, $\delta(M^\prime,M)$,
		$\delta(N^\prime,N)$ and $\delta(Y_2,Y_2^\prime)$, we have  
		\begin{align}\label{e:stable-Fredholm-index}			
			\Index(Y_1^\prime;M^\prime, N^\prime;Y_2^\prime)
			\le	\Index(Y_1;M,N;Y_2).					
		\end{align}	
		\item[(c)]
		For each tuple $(Y_1^\prime;M^\prime, N^\prime;Y_2^\prime)$ of closed linear subspaces with  $Y_1^\prime\subset M^\prime\cap N^\prime\subset M^\prime+N^\prime\subset Y_2^\prime$, and sufficiently small $\hat\delta(Y_1^\prime,Y_1)$, $\hat\delta(M^\prime,M)$,
		$\hat\delta(N^\prime,N)$ and $\hat\delta(Y_2^\prime,Y_2)$, we have  
		\begin{align}\label{e:stable-Fredholm-index-3}
			\Index(Y_1^\prime;M^\prime, N^\prime;Y_2^\prime)
			=\Index(Y_1;M,N;Y_2).	
		\end{align}
		\item[(d)] Let $m$ be an integer. Assume that $\dim (M\cap N)/Y_1$ is infinite and $\dim Y_2/(M+N)$ is finite. Then each tuple $(Y_1^\prime;M^\prime, N^\prime;Y_2^\prime)$ of closed linear subspaces with  $Y_1^\prime\subset M^\prime\cap N^\prime\subset M^\prime+N^\prime\subset Y_2^\prime$, and sufficiently small $\delta(Y_1^\prime,Y_1)$, $\delta(M,M^\prime)$,
		$\delta(N,N^\prime)$ and $\delta(Y_2^\prime,Y_2)$, we have
		\begin{align}\label{e:stable-right-semi-Fredholm}
			\Index(Y_1^\prime;M^\prime, N^\prime;Y_2^\prime)\ge m.
		\end{align}	
		\item[(e)] Let $m$ be an integer. Assume that $\dim (M\cap N)/Y_1$ is finite and $\dim Y_2/(M+N)$ is infinite. Then each tuple $(Y_1^\prime;M^\prime, N^\prime;Y_2^\prime)$ of closed linear subspaces with  $Y_1^\prime\subset M^\prime\cap N^\prime\subset M^\prime+N^\prime\subset Y_2^\prime$, and sufficiently small $\delta(Y_1,Y_1^\prime)$, $\delta(M^\prime,M)$,
		$\delta(N^\prime,N)$ and $\delta(Y_2,Y_2^\prime)$, we have  
		\begin{align}\label{e:stable-left-semi-Fredholm}
			\Index(Y_1^\prime;M^\prime, N^\prime;Y_2^\prime)\le m.		
		\end{align}  
	\end{itemize}      
\end{theorem}

\begin{rem}\label{r:stable-Fredholm-index} (a)
	By \cite[Sectioin I.4.6]{Ka95}, Theorem \ref{t:stable-Fredholm-index} follows from \cite[Theorem IV.4.30]{Ka95} when $Y_1$ is  complemented in $Y_2$ and $Y_2$ is complemented.
	\newline (b) The statements (d), (e) of Theorem \ref{t:stable-Fredholm-index} is weaker than \cite[Theorem IV.4.30]{Ka95}. The corresponding stronger results do not hold in general.
\end{rem}

Theorems \ref{t:stable-semi-Fredholm} and \ref{t:stable-Fredholm-index} are claimed to be true even under a weaker topology than gap (\cite[Theorems 4.2.1-4.2.4]{Bur08} and \cite[Theorem 2.5.1]{Bur09}). However, we have not found any published version of his paper. Moreover, our method relies on the analysis on the gaps the minimal gaps only.
So our method is easier to understand.

Though in general the stability for the indices of semi-Fredholm tuples does not hold, we still have the following

\begin{theorem}\label{t:stable-Fredholm-index-family}
	Let $X$ be a Banach space with a continuous family of semi-Fredholm tuple $\{(Y_1(s);M(s),N(s);Y_2(s))\}_{s\in J}$ of closed linear subspaces, where $J$ is a connected topological space. Then the function $s\mapsto\Index(Y_1(s);M(s),N(s);Y_2(s))$ is a constant function on $s\in J$. 	
\end{theorem}

The last problem in this paper is the study of perturbed augmented Morse index.

\begin{theorem}\label{t:perturbed-augmented-Morse-index}
	Let $X$ be a Banach space with two bounded symmetric pairs $(Q,V)$ and $(R,W)$. Let $V_0\subset V^Q$ and $W_0\subset W^R$ be two closed linear subspaces of $X$. Let $h\in\{1,-1\}$ and $c\ge 0$ be two real numbers.	
	Let $\alpha$ be a closed positive semi-definite linear subspace of $V$ with respect to $hQ$. Let $\beta$ be a closed negative semi-definite linear subspace of $V$ with respect to $hQ$ with  $V_0=\beta^{Q|\beta}$ and $\gamma(Q|_\beta)>0$. 
	Assume that there is a $Q$-orthogonal direct sum decomposition
	\begin{align}\label{e:V-a-b}
		V=\alpha\oplus\beta.
	\end{align}
	Then for each $(R,W)$ with sufficiently small $\hat\delta(V,W)$, $\delta(V_0,W_0)$ and $\delta_c(Q,R)$, we have
	\begin{align}\label{e:perturbed-augmented-Morse-index}
		m^{+}(hR)+\dim W^R/W_0\le \dim\alpha.
	\end{align} 
\end{theorem}

The following inequality for perturbed Morse index was proved by H. Li and the second author. We state it here for completeness.

\begin{proposition}\label{p:perturbed-Morse-index}\cite[Proposition 2.21]{LhcZ24}
	Let $X$ be a Banach space with two bounded symmetric pairs $(Q,V)$ and $(R,W)$. Let $c\ge 0$ be a real number.
	Let $\alpha$ be a linear subspace of $V$ with $\dim\alpha\in[1,+\infty)$, $h$ be in $\{1,-1\}$ and $h Q|_\alpha$ be positive definite. Then for each $(R,W)$ with sufficiently small $\delta(V,W)$ and $\delta_c(Q,R)$, we have
	\begin{align}\label{e:perturbed-Morse-index}
		m^+(hR)\ge \dim\alpha.
	\end{align} 	    
\end{proposition}

The proof of Theorem \ref{t:perturbed-augmented-Morse-index} is 
surprisingly nontrivial. We firstly estimate the gap between annihilators (see Lemma \ref{l:annihilator-close} below). Then we study the case of $Q|_{\alpha}=0$ (see Proposition \ref{p:perturbed-Morse-index-definite} below). Finally we are able to prove Theorem \ref{t:perturbed-augmented-Morse-index}.

Throughout this paper, we denote by $\KK$ the field of real numbers or complex numbers. We denote by $\N$, $\Z$, $\R$ and $\C$ the sets of all natural, integral, real and complex numbers respectively. By $S^1$ we denote the unit circle in the complex plane. Denote by $I_X$ the identity map on a set $X$. If there is no confusion, we will omit the subindex $X$. For a Banach space $X$, we denote by $\Bb(X)$ the set of bounded linear operators. We equip the the set of closed linear subspaces $\Ss(X)$ with {\em gap distance} $\hat\delta$ (\cite[Section IV.2.1]{Ka95}). We denote by $\gamma(M,N)$ and $\hat\gamma(M,N)$ the quantities defined by Definition \ref{d:closed-distance} below (see \cite[(IV.4.4), (IV.4.5)]{Ka95}). We denote the Fredholm index for two linear subspaces $(M,N)$ of $X$ by $\Index(M,N;X)$. We denote by $\Index(M,N)=\Index(M,N;X)$ if there is no confusion. For a symmetric form  $Q$ over $\KK$, we denote by $m^{\pm}(Q)$ and $m^0(Q)$ the Morse positive (or negative) index and the nullity of $Q$ respectively. 

The paper is organized as follows. In Section 1, we explain why we make the research, state the main results of the paper, and explain the difficulties in the proof. In Section 2, firstly we introduce some basic notions in Banach spaces. We introduce the Riesz-Schauder theory for compact operators with closed domain in a normed vector space. Then we introduce the notion of semi-compact perturbations and study the properties of them. After that, we prove Theorem \ref{t:semi-compact-perturbation} (structure theorem), Theorem \ref{t:relative-dimension-dual} (duality formula), Theorem \ref{t:index-compact-perturbations} (index formula for compact perturbation), Theorem \ref{t:stable-semi-compact-perturbations-complemented}, Theorem \ref{t:stable-relative-dimention} and Theorem \ref{t:stable-relative-dimention-family}. In Section 3, we introduce the notion of gaps and Fredholm tuples. Then we prove the splitting theorem (Theorem \ref{t:existence-finite-dimension-embeding-close-oposite} below), Theorem \ref{t:stable-semi-Fredholm}, Theorem \ref{t:stable-Fredholm-index} and Theorem \ref{t:stable-Fredholm-index-family}.  In Section 4, we study the perturbed augmented Morse index and prove Theorem \ref{t:perturbed-augmented-Morse-index}. 

We would like to thank the referees of this paper for their critical reading and
very helpful comments and suggestions.

\section{Stability theorems for semi-compact perturbations}\label{s:stable-semi-compact- perturbations}

\subsection{Fredhom pairs of closed linear subspaces}\label{ss:Fredholm}

Let $X$ be a Banach space. We denote by $\Ss(X)$ the set of closed linear subspaces of $X$ and $\Ss^c(X)$ the set of complemented closed linear subspaces of $X$ respectively. Let $M,N$ be two closed linear subspaces of $X$, i.e., $M,N\in\Ss(X)$.
Denote by $\romS_M$ the unit
sphere of $M$. We recall three common definitions of distances in $\Ss(X)$ (see also
\cite[Sections IV.2.1 and IV.4.1]{Ka95}):
\begin{itemize}
	
	\item  the \textit{Hausdorff metric} $\hat d$;
	
	\item the \textit{aperture (gap distance)} $\hat\delta$, that is not a metric since it does not in general
	satisfy the triangle inequality, but defines the same topology as the metric $\hat d$, called
	\textit{gap topology}, and is easier to estimate than $\hat d$; and
	
	\item the \textit{angular distance (minimum gap)} $\wh\gamma$, that is useful in our estimates, though not
	defining any suitable topology.
	
\end{itemize}

\begin{definition}[The gap between linear subspaces]\label{d:closed-distance}	
	(a) We set \begin{multline*}
		\hat d(M,N)\ =\
		d(\romS_M,\romS_N)\\
		:= \begin{cases} \max\left\{\begin{matrix}\sup\limits_{u\in \romS_M}\dist(u,\romS_N),\\ \sup\limits_{u\in
					\romS_N}\dist(u,\romS_M)\end{matrix}\right\},&\text{ if both $M\ne 0$ and $N\ne 0$},\\
			0,&\text{ if $M=N=0$},\\
			2,&\text{ if either $M= 0$ and $N\ne 0$ or vice versa}.
		\end{cases}
	\end{multline*}
	\newline (b) We set
	\begin{align*}
		\delta(M,N)\ :=&\ \begin{cases}\sup\limits_{u\in \romS_M}\dist(u,N),& \text{if $M\ne\{0\}$},\\
			0,& \text{if $M=\{0\}$},\end{cases}\\
		=&\ \delta(\bar M,\bar N),\\
		\hat\delta(M,N)\ :=&\ \max\{\delta(M,N),\delta(N,M)\}=\ \hat\delta(\bar M,\bar N).
	\end{align*}
	$\hat\delta(M,N)$ is called the {\em gap} between $M$ and $N$.
	\newline (c) Assume that $M\cap N$ is closed. We set
	\begin{align*}
		\gamma(M,N)\ :=&\ \begin{cases} \inf\limits_{u\in M\setminus N}\frac{\dist(u,N)}{\dist(u,M\cap N)}\ (\le
			1), & \text{if $M\not\subset N$},\\
			1, & \text{if $M\subset N$},\end{cases}\\
		\hat\gamma(M,N)\ :=&\ \min\{\gamma(M,N),\gamma(N,M)\}.
	\end{align*}
	$\hat\gamma(M,N)$ is called the {\em minimum gap} between $M$ and $N$. If $M\cap N=\{0\}$, we have
	\[\gamma(M,N)\ =\ \inf\limits_{u\in \romS_M}\dist(u,N).\]
\end{definition}

In this paper we shall impose the gap topology on the space $\Ss(X)$ of all closed linear
subspaces of a Banach space $X$ and its subset $\Ss^c(X)$ of complemented subspaces.

We recall the following two results on finite-dimensional variation. For the second see
\cite[Proposition 11.4]{Brezis:2011}. For the first see (\cite[Lemma III.1.9]{Ka95}.

\begin{proposition}[Finite extension]\label{p:finite-extension} Let $X$ be a Banach space and $M$ be a closed subspace of $X$.
	Let $M^{\prime}\supset M$ be a linear subspace of $X$ with $\dim M^{\prime}/M<+\infty$. Then we
	have
	\newline (a) $M^{\prime}$ is closed, and
	\newline (b) $M'\in\Ss^c(X)$ if and only if $M\in\Ss^c(X)$.
\end{proposition}

\begin{definition}\label{d:fredholm-pair}
	(a) The space of (algebraic) \emph{Fredholm pairs} of linear subspaces of a vector space $X$ is
	defined by
	\begin{equation}\label{e:fp-alg}
		\Ff^{2}_{\operatorname{alg}}(X):=\{(M,N);  \dim (M\cap N)  <+\infty\;\; \text{and} \dim
		X/(M+N)<+\infty\}
	\end{equation}
	with
	\begin{equation}\label{e:fp-index}
		\Index(M,N;X)=\Index(M,N):=\dim(M\cap N) - \dim X/(M+N).
	\end{equation}
	
	\noi (b) In a Banach space $X$, the space of (topological) \emph{Fredholm pairs} is defined by
	\begin{multline}\label{e:fp}
		\Ff^{2}(X):=\{(M,N)\in\Ff^2_{\operatorname{alg}}(X); M,N, \tand M+N \subset X \text{ closed}\}.
	\end{multline}
	A pair $(M,N)$ of closed subspaces is called {\em semi-Fredholm} if $M+N$ is closed, and at least
	one of $\dim (M\cap N)$ and $\dim X/(M+N)$ is finite.
	
	\noi (c) Let $X$ be a Banach space, $M\in\Ss(X)$ and $k\in\Z$. We define
	\begin{align}
		\label{e:fp-M}\Ff_M(X):&=\{N\in\Ss(X);(M,N)\in\Ff^2(X)\},\\
		\label{e:fp-kM}\Ff_{k,M}(X):=&\{N\in\Ss(X);(M,N)\in\Ff^2(X),\Index(M,N)=k\}.
	\end{align}
\end{definition}

\begin{rem}\label{r:redholm-pairs}
	Actually, in Banach space the closeness of $\la+\mu$ follows from its finite codimension in $X$ in
	combination with the closeness of $\la,\mu$ (see \cite[Remark A.1]{BoFu99} and \cite[Problem
	IV.4.7]{Ka95}).
\end{rem}

As in \cite[Sectioin III.1.4]{Ka95}, the {\em adjoint space} $X^*$ of a Banach space is defined to be the set of bounded semilinear forms on $X$. Let $M$ be a subset of the Banach space $X$. We define
\begin{align}\label{e:M-perp}
	M^\perp\ :=\{f\in X^*;\;f(x)=0\text{ for all } x\in M\}.
\end{align}

We have the following (cf.\cite[Problem IV.4.6]{Ka95}).

\begin{lemma}\label{l:finite-diff-index} Let $X$ be a vector space and $M^{\prime},M,N$ be linear subspaces. Assume that
	$M^{\prime}\supset M$ holds.
	\begin{itemize}
		\item[(a)] If $\dim(M^\prime\cap N)<+\infty$, we have $\dim(M\cap N)<+\infty$ and $\Index(M,N)=\Index(M^{\prime},N)-\dim M^\prime/M$.
		\item[(b)] If $\dim X/(M+N)<+\infty$, we have $\dim X/(M^\prime+N)<+\infty$ and $\Index(M^{\prime},N)=\Index(M,N)+\dim M^\prime/M$.
	\end{itemize} 
\end{lemma}

\begin{proof}
	We have 
	\begin{align*}
		\frac{X/(M+N)}{(M^\prime+N)/(M+N)}
		\simeq&X/(M^\prime+N),\quad\tand\\
		\frac{M^\prime/M}{(M^\prime\cap N)/(M\cap N)}
		\simeq&\frac{M^\prime/M}{(M+M^\prime\cap N)/M}
		\simeq \frac{M^\prime}{M+M^\prime\cap N}.
	\end{align*}
	By Lemma \cite[A.1.1]{BoZh18} we have
	\begin{align*}
		\frac{M^\prime+N}{M+N}
		\simeq\frac{M^\prime}{M^\prime\cap(M+N)}=\frac{M^\prime}{M+M^\prime\cap N}.
	\end{align*}
	Then (a), (b) follows.
\end{proof}

\begin{lemma}\label{l:subspace-index} 
	Let $X$ be a Banach space with a left Fredholm pair $(M,N)$ of closed linear subspaces. Let
	$M^{\prime}\subset M$ and $N^{\prime}\subset N$ be closed linear subspaces. Then $(M^\prime,N^\prime)$ is left Fredholm, and we have
	\begin{align}\label{e:subspace-index}
		\Index(M^{\prime},N^\prime)=\Index(M,N)-\dim M/M^\prime-\dim N/N^\prime.
	\end{align}
\end{lemma}
\begin{flushleft}
	\textit{}
\end{flushleft}
\begin{proof}
	By \cite[Theorem IV.4.2]{Ka95}, the lemma is clear if $M\cap N=\{0\}$.
	
	In the general case, we have $\dim(M^{\prime}\cap N^\prime)\le \dim(M\cap N)<+\infty$. By introducing the quotient space $X/(M\cap N)$, the space $(M^\prime+N^\prime+M\cap N)/(M\cap N)$ is closed in $X/(M\cap N)$. Then the space $M^\prime+N^\prime+M\cap N$ is closed in $X$. By \cite[Problem IV.4.7]{Ka95}, the space $M^\prime+N^\prime$ is closed in $M^\prime+N^\prime+M\cap N$. Then the space $M^\prime+N^\prime$ is closed in $X$. Thus $(M^\prime,N^\prime)$ is left Fredholm. By Lemma \ref{l:finite-diff-index} we obtain \eqref{e:subspace-index}.	
\end{proof}

\subsection{Riesz-Schauder theory}
Let $X$ be a normed vector space. Let $K\colon X\supset\dom(K)\to X$ be a compact linear operator with closed domain. In this subsection we describe Riesz-Schauder theory for the spectrum of $K$.

\begin{lemma}\label{l:compact-spectrum-C}
	Let $X$ be a Banach space. Let $K\colon X\supset\dom(K)\to X$ be a compact linear operator such that $\dom(K)$ is a closed subset of $X$. Then for each $\lambda\ne\{0\}$, we have 
	\begin{align}\label{e:compact-index}
		\Index(\lambda I-K)=-\dim X/\dom(K).
	\end{align}
	Moreover, the spectrum of $K$ is $\C$ if $\dom(K)$ is a proper closed subset of $X$.
\end{lemma} 

\begin{proof}
	By Theorem \cite[Theorem IV.5.26]{Ka95} we have
	\begin{align*}
		\Index(\lambda I-K)
		=\Index(\lambda I_{\dom(K)})=-\dim X/\dom(K).
	\end{align*}
	Then we have
	\begin{align}\label{e:compact-range-small}
		\dim X/\image(\lambda I-K)\ge\dim X/\dom(K).
	\end{align}
	If $\dom(K)$ is a proper closed subset of $X$, the spectrum of $K$ contains $\lambda$ and $0$. Thus the spectrum of $K$ is $\C$.
\end{proof}

By Lemma \ref{l:compact-spectrum-C}, the spectrum of $K$ is always $\C$ if $K$ is a compact linear operator of a Banach space $X$ with proper closed domain. In the following we study the eigenvalues of $K$.

We need a lemma.

\begin{lemma}\label{l:compact-compact}
	Let $X$ be a normed vector space. Let $K\colon X\supset\dom(K)\to X$ be a compact linear operator. Let $n\ge 0$ be an integer. Then the operator $I-(I-K)^n$ is compact.
\end{lemma}

\begin{proof}
	Since $(I-K)^{n+1}=(I-K)^n-K(I-K)^n$, by induction the operator $I-(I-K)^n$ is compact.
\end{proof}

\begin{proposition}\label{p:Fredholm-compact}
	Let $X$ be a normed vector space. Let $K\colon X\supset\dom(K)\to X$ be a compact linear operator with closed domain. Then we have the following.
	\begin{itemize}
		\item [(a)] There exists a constant $C>0$ such that for each $x\in\dom(K)$, there holds
		\begin{align}\label{e：compact-dist-kernel}
			\dist(x,\ker(I-K))\le C\|(I-K)x\|.
		\end{align}
		\item[(b)] The range $\image(I-K)$ is closed.
		\item[(c)] There exists an integer $n\ge 0$ such that the operator $I-K\colon\image(I-K)^n\to\image(I-K)^n$ is a bounded operator with bounded inverse. Moreover, we have 
		\begin{align}\label{e:compact-kernel-image}
			\dom(I-K)^n=&\ker(I-K)^n\oplus \image(I-K)^n,\\
			\label{e:compact-eigenspace}
			\ker(\lambda I-K)^n=&\ker(\lambda I-K)^{n+1}.
		\end{align}	
	\end{itemize}
\end{proposition}

\begin{proof} 
	(a), (b) The proof is the same as that of \cite[Theorem 5.3]{GiTr01} and we omit it here.
	\newline (c) By Step 3 in the proof of \cite[Theorem 5.3]{GiTr01}, there exists an integer $n\ge 0$ such that the operator $I-K\colon\image(I-K)^n\to\image(I-K)^n$ is a bounded surjective operator. By \cite[Theorem 5.3]{GiTr01}, the operator $I-K\colon\image(I-K)^n\to\image(I-K)^n$ is a bounded operator with bounded inverse. Then we have $\ker(I-K)^n\cap \image(I-K)^n=\{0\}$. Then \eqref{e:compact-kernel-image} follows from Step 3 in the proof of \cite[Theorem 5.3]{GiTr01}.
	Then we have 
	\begin{align*}
		\ker(\lambda I-K)^n=&\ker((\lambda I-K)^n\colon\dom(\lambda I-K)^n\to\image(\lambda I-K)^n)\\
		=&\ker(\lambda I-K)^{2n}=\ker(\lambda I-K)^{n+1}.
	\end{align*}
\end{proof}

\begin{proposition}\label{p:eigenvalues-compact}
	Let $X$ be a normed vector space. Let $K\colon X\supset\dom(K)\to X$ be a compact linear operator. Then the set of eigenvalues of $K$ is a countable set with at most one limited point $0$. Each nonzero eigenvalue of $K$ has finite geometric multiplicity. Each nonzero eigenvalue of $K$ has finite algebraic multiplicity if $\dom(K)$ is closed.
\end{proposition}

\begin{proof} 
	The proof is the same as that of \cite[Theorem 5.5]{GiTr01} except the last statement.
	
	Let $\lambda$ be a nonzero eigenvalue of $K$. Then the proof of \cite[Theorem 5.5]{GiTr01} shows that $\dim\ker(\lambda I-K)<+\infty$.
	By Lemma \ref{l:compact-compact}, the operator $I-(I-K)^n$ is compact.
	By Proposition \ref{p:Fredholm-compact}.c, there exists an integer $n\ge 0$ such that such that the operator $\lambda I-K\colon\image(I-K)^n\to\image(\lambda I-K)^n$ is a bounded operator with bounded inverse and there holds
	\begin{align*}
		\dom(\lambda I-K)^n=\ker(\lambda I-K)^n\oplus \image(\lambda I-K)^n.
	\end{align*}

	By Lemma \ref{l:compact-compact}, $\lambda^n I-(\lambda I-K)^n$ is compact.
    By the proof of \cite[Theorem 5.5]{GiTr01}, we have
    \begin{align*}
    	\dim\ker((\lambda I-K)^n\colon\dom(\lambda I-K)^n\to\image(\lambda I-K)^n)<+\infty.
    \end{align*}
    By Proposition \ref{p:Fredholm-compact}.c, $\lambda$ is an eigenvalue of $K$ with finite algebraic multiplicity.
\end{proof}

\subsection{Relative dimension between closed linear subspaces}\label{ss:relative-dimension}
Firstly we define the relative dimension.

\begin{definition}\label{d:semi-compact-perturtation} Let $X$ be a Banach space 
with two closed linear subspaces $\lambda$, $\mu$.
\begin{itemize}
	\item[(a)] We define $\lambda\sim^{sc}\mu$ if there is a compact linear operator $K\in\Bb(\lambda,X)$ such that $(I+K)\lambda\subset\mu$, and call $\mu$ a {\em right semi-compact perturbation} of $\lambda$. In this case, we define the {\em relative dimension} $[\lambda-\mu]\in\Z\cup\{-\infty\}$ by
	\begin{align}\label{e:semi-compact-perturtation} 
		[\lambda-\mu]\ :=[\lambda-\mu;K]\ :=\Index(I+K\colon\lambda\to\mu).
	\end{align}
    We define $[\mu-\lambda]=+\infty$ if $[\lambda-\mu]=-\infty$.
    \item[(b)] We define $\lambda\sim^{gsc}\mu$ if there is a compact linear operator $K\in\Bb(X)$ such that $(I+K)\lambda\subset\mu$, and call $\mu$ a {\em right global semi-compact perturbation} of $\lambda$.
    \item[(c)] We define $\lambda\overset{sc}{\sim}\mu$ if either $\lambda\sim^{sc}\mu$ or $\mu\sim^{sc}\lambda$, and call $\mu$ a {\em semi-compact perturbation} of $\lambda$. We define $\lambda\overset{gsc}{\sim}\mu$ if either $\lambda\sim^{gsc}\mu$ or $\mu\sim^{gsc}\lambda$, and call $\mu$ a {\em global semi-compact perturbation} of $\lambda$.
    \item[(d)] We define $\lambda\sim^c\mu$ if there is a compact linear operator $K\in\Bb(\lambda,X)$ such that $(I+K)\lambda\subset\mu$ and the operator $I+K\colon\lambda\to\mu$ is Fredholm, and call $\lambda$ a {\em compact perturbation} of $\mu$.
    We define $\lambda\sim^{gc}\mu$ if there is a compact linear operator $K\in\Bb(X)$ such that $(I+K)\lambda\subset\mu$ and the operator $I+K\colon\lambda\to\mu$ is Fredholm, and call $\lambda$ a {\em global compact perturbation} of $\mu$.
\end{itemize} 
\end{definition}

The relative dimension is well-defined. More precisely, we have

\begin{lemma}\label{l:semi-compact-perturtation}
	Let $X$ be a Banach space
	with three closed linear subspaces $\alpha$, $\beta$ and $\gamma$. 
	\begin{itemize}
		\item[(a)] Assume that $\alpha\sim^{sc}\beta$ holds. Then the relative dimension $[\alpha-\beta]$ is well-defined.
		\item[(b)] Assume that $\alpha\sim^{sc}\beta$ and $\beta\sim^{sc}\gamma$ hold. Then  we have $\alpha\sim^{sc}\gamma$, and
		\begin{align}\label{e:semi-compact-perturtation-abc}
			[\alpha-\gamma]=[\alpha-\beta]+[\beta-\gamma].
		\end{align} 
		\item[(c)] Assume that $\alpha\sim^{gsc}\beta$ and $\beta\sim^{gsc}\gamma$ hold. Then we have $\alpha\sim^{gsc}\gamma$.
		\item[(d)] For each compact operator $K\in\Bb(\alpha,X)$ with $\|K\|<1$, we have $[\alpha-(I+K)\alpha]=0$.
	\end{itemize}  
\end{lemma}

\begin{proof}
	(a) Let $K_1,K_2\in\Bb(\lambda,X)$ be two compact operator with $(I+K_i)\alpha\in\beta$ for $i=1,2$. 
	Since $I|_\lambda\in\Bb(\lambda,X)$ is left Fredholm (i.e., Fredholm with finite dimensional kernel), by \cite[Theorem IV.5.26]{Ka95}, $I+K_1,I+K_2\in\Bb(\lambda,X)$ are left Fredholm. Then $I+K_1,I+K_2\in\Bb(\lambda,\mu)$ are left Fredholm.
	Since $K_1-K_2\colon\alpha\to\beta$ is compact, by \cite[Theorem IV.5.26]{Ka95}, we have 
	\begin{align*}
		[\alpha-\beta;K_1]=[\alpha-\beta;K_2].
	\end{align*}
    \newline (b) Since $\alpha\sim^{sc}\beta$ and $\beta\sim^{sc}\gamma$ hold, there are compact linear operators $K\in\Bb(\lambda,X)$ and $L\in\Bb(\mu,X)$ such that 
    \begin{align*}
    	(I+K)\alpha\subset\beta\text{ and }(I+L)\beta\subset\gamma.
    \end{align*}
    Then we have $(I+L)(I+K)\alpha\subset\gamma$. Since $(I+L)(I+K)-I=L(I+K)+K\in\Bb(\lambda,X)$ is compact, we have $\alpha\sim^{sc}\gamma$. By \cite[Exercises B.11]{GWh78} we have
    \begin{align*}
    	[\alpha-\gamma]&=[\alpha-\gamma;(I+L)(I+K)-I]\\
    	&=[\alpha-\beta;K]+[\beta-\gamma;L]\\
    	&=[\alpha-\beta]+[\beta-\gamma].
    \end{align*}    
    \newline (c) Since $\alpha\sim^{gsc}\beta$ and $\beta\sim^{gsc}\gamma$ hold, there are compact linear operators $K,L\in\Bb(X)$ such that 
    \begin{align*}
    	(I+K)\alpha\subset\beta\text{ and }(I+L)\beta\subset\gamma.
    \end{align*}
    Then we have $(I+K)(I+L)\alpha\subset\gamma$. Since $(I+K)(I+L)-I\in\Bb(X)$ is compact, we have $\alpha\sim^{gsc}\gamma$.
    \newline (d) Since $\|K\|<1$, we have $\ker(I+K)=0$. Then we have
    \begin{align*}
    	[\alpha-(I+K)\alpha]=\dim\ker(I+K)=0.
    \end{align*}
\end{proof}

The relative dimension is preserved by a bounded linear embedding.

\begin{lemma}\label{l:relative-dimension-embedding}
	Let $X$ be a Banach space 
	with two closed linear subspaces $\lambda$, $\mu$. Assume that $\lambda\sim^{sc}\mu$. Let $Y$ be a Banach space and $j\in\Bb(X,Y)$ be such that $jX$ is closed in $Y$ and $j^{-1}\in\Bb(jX,X)$. Then we have $j\lambda\sim^{sc}j\mu$, and
	\begin{align}\label{le:relative-dimension-isomorphism}
		[j\lambda-j\mu]=[\lambda-\mu].
	\end{align}
\end{lemma}

\begin{proof}
	Since $\lambda\sim^{sc}\mu$, there is a compact operator $K\in\Bb(\lambda,X)$ such that $(I_X+K)\lambda\subset\mu$. Then $jKj^{-1}\in\Bb(j\lambda,X)$ is compact, $I_Y+jKAj^{-1}=j(I_X+K)j^{-1}$, and $(I+jKj^{-1})j\lambda\subset j\mu$. Then we have $j\lambda\sim^{sc}j\mu$, and
	\begin{align*}
		[j\lambda-j\mu]=&\Index\left((I_Y+jKj^{-1})\colon j\lambda\to j\mu\right)\\
		=&\Index(I_X+K\colon \lambda\to \mu)=
		[\lambda-\mu].
	\end{align*}
	
\end{proof}

The notion of global compact perturbation coincides with that of compact perturbation defined by \cite[Definition A.7.3(b)]{BoZh18}. We begin with some preparations.

\begin{lemma}\label{l:Fredholm-closed}
	Let $X$ and $Y$ be two Banach spaces. Let $F\in\Bb(X,Y)$ be a bounded operator with
	closed range. Then for each linear subspace $\lambda$ of $X$ such that $\lambda+\ker F$ is closed in $X$, the space $F\lambda$ is a closed linear subspace of $Y$. Moreover, we have linear isomorphisms
	\begin{align*}
		(F\lambda)/(F\lambda_1)&\simeq(\lambda+\ker F)/(\lambda_1+\ker F)\\
		&\simeq(\lambda/(\lambda\cap\ker F))/(\lambda_1/(\lambda_1\cap\ker F))		
	\end{align*}
    for each pair $(\lambda,\lambda_1)$ of linear subspaces of $X$ with $\lambda_1\subset\lambda$.
\end{lemma}

\begin{proof}
	1. The operator $F$ induce a linear isomorphism 
	\begin{align*}
		\tilde F\in\Bb(X/\ker F,\ran F).
	\end{align*} 
    Since $\ran F$ is closed in $Y$, the space $F\lambda=\tilde F(\lambda+\ker F)$ is a closed linear subspace of $\ran F$.  Then the space $F\lambda$ is a closed linear subspace of $Y$.
    \newline 2. Since $\tilde F$ is a linear isomorphism, we have
    \begin{align*}
    	(F\lambda)/(F\lambda_1)
    	&\simeq (\tilde F(\lambda+\ker F))/(\tilde F(\lambda_1+\ker F))\\
    	&\simeq(\lambda+\ker F)/(\lambda_1+\ker F)\\
    	&\simeq((\lambda+\ker F)/\ker F)/((\lambda_1+\ker F)/\ker F)\\
    	&\simeq(\lambda/(\lambda\cap\ker F))/(\lambda_1/(\lambda_1\cap\ker F)).		
    \end{align*}
\end{proof}

\begin{corollary}\label{c:compact-semi-Fredholm}
	Let $X$ be a Banach space 
	with two closed linear subspaces $\lambda$, $\mu$.
	Let $K\in\Bb(X)$ a compact linear operator such that $(I+K)\lambda\subset\mu$. Then we have $\dim\ker(I+K)|_\lambda\le\dim\ker(I+K)$ and $(I+K)\lambda$ is closed.
\end{corollary}

\begin{lemma}\label{l:complemented-finite}
	Let $X$ be a Banach space
	with two closed linear subspaces $\lambda$ and $V$. Assume that $\dim V<+\infty$ and $V\cap\lambda=\{0\}$. Then there is a closed linear subspace $\mu\supset\lambda$ of $X$ such that $X=V\oplus\mu$.
\end{lemma}

\begin{proof}
	Let $v_1,\ldots,v_n$ be a linear basis of $V$. We define 
	\begin{align*}
		f_k:\lambda+ \C v_k\to\C
	\end{align*}
    by $f_k(v_k)=1$ and $f_k|_\lambda=0$ for each $k=1,\ldots,n$. By Hahn-Banach theorem, there is an extension $\tilde f_k\in X^*$ of $f_k$. 
	Set $\mu\ :=\cap_{k=1}^n\ker\tilde f_k$.
	Then $\mu$ is a closed linear subspace of $X$ with $\lambda\subset\mu$ and $X=V\oplus\mu$.
\end{proof}

\begin{definition}\label{d:space-finite-change}
	Let $X$ be a linear space 
	with two linear subspaces $\lambda$, $\mu$. We define $\lambda\sim^f\mu$ if $\dim\lambda/(\lambda\cap\mu)<+\infty$. We call $\mu$ a {\em finite change} of $\lambda$ (see \cite[p.273]{Neu68}) if $\lambda\sim^f\mu$ and $\mu\sim^f\lambda$. We define $\lambda\overset{sf}{\sim}\mu$ if either $\dim\lambda/(\lambda\cap\mu)<+\infty$ or $\dim\mu/(\lambda\cap\mu)<+\infty$, and call $\mu$ a {\em semi-finite change} of $\lambda$. In this case, we define
	the {\em relative dimension} $[\lambda-\mu]$ by
	\begin{align*}
		[\lambda-\mu]\ :=\dim\lambda/(\lambda\cap\mu)-\dim\mu/(\lambda\cap\mu).
	\end{align*}
\end{definition}

We have the following characterization of semi-finite change.

\begin{lemma}\label{l:finite-change}
	Let $X$ be a Banach space 
	with two closed linear subspaces $\lambda$, $\mu$. Then 
	there holds $\lambda\sim^f\mu$ if and only if there is a finite rank operator $F\in\Bb(X)$ such that $(I+F)\lambda\subset\mu$, if and only if there is a finite rank operator $F\in\Bb(\lambda,X)$ such that $(I+F)\lambda\subset\mu$. In this case we have $\lambda\sim^{gc}\mu$.
\end{lemma}

\begin{proof}
	1. {\em Necessity}. 
	
	Since $\dim(\lambda/\lambda\cap\mu)<+\infty$ holds, there is a finite linear subspace $V$ of $\lambda$ such that 
	\begin{align*}
		\lambda=V\oplus\lambda\cap\mu.
	\end{align*}
	By Lemma \ref{l:complemented-finite} there is a closed linear subspace $\alpha$ of $X$ containing $\lambda\cap\mu$ such that
	\begin{align*}
		X=V\oplus\alpha.
	\end{align*}
	Let $P$ be the projection of $X$ onto $V$ along $\alpha$. Then $F\colon\ =-P\in\Bb(X)$ is a finite rank operator and we have 
	\begin{align*}
		(I+F)\lambda=\lambda\cap\mu\subset\mu.
	\end{align*}
	Since $F$ is compact, we have $\lambda\sim^{gc}\mu$.
	\newline 2. {\em Sufficiency}.
	
	Assume that there is a finite rank operator $F\in\Bb(\lambda,X)$ such that $(I+F)\lambda\subset\mu$. Then we have
	\begin{align*}
		\lambda\cap\mu\supset\lambda\cap((I+F)\lambda)\supset\lambda\cap\ker F.
	\end{align*}
	Then we have
	\begin{align*}
		\dim\lambda/(\lambda\cap\mu)\le&\dim\lambda/(\lambda\cap\ker F)\\
		=&\dim(\lambda+\ker F)/\ker F<+\infty.
	\end{align*}	
\end{proof}

We have the following structure theorem for semi-compact perturbations.

\begin{theorem}\label{t:semi-compact-perturbation}
	Let $X$ be a Banach space
	with two closed linear subspaces $\lambda$ and $\mu$. Then we have the following.
	\begin{itemize}
		\item[(a)] There holds $\lambda\sim^{sc}\mu$ if and only if there are closed subspaces $\lambda_1$ of $\lambda$, $\mu_1$ of $\mu$ respectively and a compact linear operator $K_1\in\Bb(\lambda,X)$ such that \begin{align}\label{e:sc-finite-codimension}
			\dim\lambda/\lambda_1<+\infty,
			\text{ and }
			(I+K_1)\lambda_1=\mu_1.
		\end{align}        
		Moreover, $\lambda_1$, $\mu_1$ and $K_1$ can be chosen such that $K_1|_{\lambda_1}=K|_{\lambda_1}$, $I+K_1$ is an injection, and $I+sK_1\colon\lambda_1\to X$ is an injection for each $s\in[0,1]$. 
		\item[(b)] There holds $\lambda\sim^{gsc}\mu$ if and only if there are closed subspaces $\lambda_1$ of $\lambda$, $\mu_1$ of $\mu$ respectively and a compact linear operator $K_1\in\Bb(X)$ such that \eqref{e:sc-finite-codimension} holds.        
		Moreover, $\lambda_1$, $\mu_1$ and $K_1$ can be chosen such that $I+K_1\in\Bb(X)$ is a linear isomorphism and $I+sK_1\in\Bb(\lambda_1,X)$ is an injection for each $s\in[0,1]$.
	\end{itemize}  
\end{theorem}

\begin{proof}
	1. {\em Necessity in (a)}. 
	
	Since $\lambda\sim^{sc}\mu$, there is a compact linear operator $K\in\Bb(\lambda,X)$ such that $(I+K)\lambda\subset\mu$, we can choose $\lambda_1=\lambda$, $\mu_1=(I+K)\lambda$ and $K_1=K$. 
	
	We need to prove that we can choose $\lambda_1$, $\mu_1$ and $K_1\in\Bb(\lambda,X)$ such that $K_1|_{\lambda_1}=K|_{\lambda_1}$, $I+sK_1\colon\lambda_1\to X$ is an injection for each $s\in[0,1]$.
			
	Set $V\ :=\sum_{s\in[0,1]}\ker(I+sK)$. By Proposition \ref{p:eigenvalues-compact}, we have $\dim V<+\infty$. 
	By Lemma \ref{l:complemented-finite}, there is a closed linear subspace $\lambda_1$ of $\lambda$ such that
	\begin{align*}
		\lambda=V\oplus\lambda_1.
	\end{align*}
    
    Set $\mu_1\ :=(I+K)\lambda_1$. Then $I+K\colon\lambda_1\to\mu_1$ is a linear isomorphism, and we have $\dim X/\mu_1=\dim X/\lambda_1\ge\dim V$. Then there exists an injection $A\in\Bb(\lambda_1,X)$ such that $\image A\cap\mu_1=\{0\}$. Define the operator $K_1\in\Bb(\lambda,X)$ by 
    $K_1(x+y)\ :=(A-I)x+Ky$ 
    for each $x\in V$ and $y\in \lambda_1$. Then $K_1|_{\lambda_1}=K|_{\lambda_1}$, $K_1$ is compact $I+K_1$ is an injection, and $I+sK_1\colon\lambda_1\to X$ is an injection for each $s\in[0,1]$. By Lemma \ref{l:Fredholm-closed}, the space $\mu_1$ is closed. 
	\newline 2. {\em Sufficiency in (a)}. 
	
	Since $\lambda_1\sim^{sc}\mu_1$, by Step 1, there are closed subspaces $\lambda_2$ of $\lambda_1$ and $\mu_2$ of $\mu_1$ and a compact linear operator $K_2\in\Bb(\lambda_1,X)$ such that 
	\begin{align*}
		\dim\lambda_1/\lambda_2<+\infty\text{ and }
		(I+K_2)\lambda_2=\mu_2.
	\end{align*} 
    Moreover,  $I+K_2\colon\lambda_2\to\mu_2$ is a linear isomorphism. 
    
    Since $\dim \lambda/\lambda_2<+\infty$, there is a finite dimensional linear subspace $W$ of $\lambda$ such that 
    \begin{align*}
    	\lambda=W\oplus\lambda_2.
    \end{align*}    
    Define $K\in\Bb(\lambda,X)$ by $K(x+y)\ :=-x+K_2y$ for $x\in W$ and $y\in \lambda_2$. Then $K$ is compact and $(I+K)\lambda=\mu_2\subset\mu$. Thus $\lambda\sim^{sc}\mu$ holds. 
    \newline 3. {\em Necessity in (b)}. 
    
    Since $\lambda\sim^{gsc}\mu$, there is a compact linear operator $K\in\Bb(X)$ such that $(I+K)\lambda\subset\mu$, we can choose $\lambda_1=\lambda$, $\mu_1=(I+K)\lambda$ and $K_1=K$. 
    
    We need to prove that we can choose $\lambda_1$, $\mu_1$ and $K_1$ such that $I+K_1\in\Bb(X)$ has bounded inverse and $I+sK_1\colon\lambda_1\to X$ is an injection for each $s\in[0,1]$. 
    
    Set $V_1\ :=\sum_{s\in[0,1]}\ker(I+sK)$. Since $K$ is compact, by Proposition \ref{p:eigenvalues-compact} we have $\dim V_1<+\infty$. 
    
    By Lemma \ref{l:complemented-finite}, there are closed linear spaces $\lambda_1$ of $\lambda$ and $X_1$ of $X$ respectively such that
    \begin{align*}
    	\lambda=V_1\cap\lambda\oplus\lambda_1,\quad
    	X=V_1\oplus X_1,\text{ and } \lambda_1\subset X_1.
    \end{align*}
    By \cite[Theorem IV.5.26]{Ka95}, we have
    \begin{align*}
    	\dim X/(I+K)X_1=&-\Index(I+K\colon X_1\to X)\\
    	=&-\Index(I\colon X_1\to X)=\dim X/X_1.
    \end{align*}
    Then there is a finite dimensional linear subspace $W_1$ of $X$ such that $X=W_1\oplus\ran (I+K)X_1$. Moreover, there is a linear isomorphism 
    $A_1\colon V_1\to W_1$. 
    
    Set $\mu_1\ :=(I+K)\lambda_1$. Then $I+K\colon\lambda_1\to\mu_1$ is a linear isomorphism. Define the operator $K_1\in\Bb(X)$ by $K_1(x+y)\ :=(A_1-I)x+Ky$ for each $x\in V_1$ and $y\in X_1$. Then $K_1$ is compact, $I+K_1$ is an injection and $I+sK_1\colon\lambda_1\to X$ is an injection for each $s\in[0,1]$. By Lemma \ref{l:Fredholm-closed}, the space $\mu_1$ is closed. By \cite[Theorem IV.5.26]{Ka95}, we have $I+K_1$ is a Fredholm operator of index $0$. Then $I+K_1$ has bounded inverse.
    \newline 4. {\em Sufficiency of (b)}. 
    
    Since $\lambda_1\sim^{gsc}\mu_1$, by Step 1, there are a closed linear subspace $\lambda_2$ of $\lambda_1$, a closed linear subspace $\mu_2$ of $\mu_1$ and a compact linear operator $K_2\in\Bb(X)$ such that 
    \begin{align*}
    	\dim\lambda_1/\lambda_2<+\infty\text{ and }
    	(I+K_2)\lambda_2=\mu_2.
    \end{align*} 
    Moreover, the operators $I+K_2$ and $I+K_2\colon\lambda_2\to\mu_2$ are linear isomorphisms. 
    
    Since $K_2$ is compact, by \cite[Theorem IV.5.26]{Ka95}, $I+K_2$ is Fredholm of index $0$. There is a finite dimensional linear subspace $V_2$ of $X$ such that 
    \begin{align*}
    	\ker(I+K_2)=V_2\oplus\ker(I+K_2)\cap\lambda.
    \end{align*} 
    Then we have $V_2\cap\lambda\subset V_2\cap\ker(I+K_2)\cap\lambda=\{0\}$. Since $\ker(I+K_2)\cap\lambda_2=\{0\}$ holds, there is a finite dimensional linear subspace $V_3\supset \ker(I+K_2)\cap\lambda$ of $\lambda$ such that 
    \begin{align*}
    	\lambda=V_3\oplus\lambda_2.
    \end{align*}  
    Then we have $V_2\cap V_3\subset V_2\cap\lambda=\{0\}$ and 
    \begin{align*}
    	(V_2+V_3)\cap\lambda_2&=(V_2\cap(V_3+\lambda_2)+V_3)\cap\lambda_2\\
    	&=(V_2\cap\lambda+V_3)\cap\lambda_2\\
    	&=V_3\cap\lambda_2=\{0\}.
    \end{align*}
    By Lemma \ref{l:complemented-finite}, there is a closed linear subspaces $X_2\supset \lambda_2$ of $X$ such that \begin{align*}
    	X=V_2\oplus V_3\oplus X_2.
    \end{align*}
    Define $K\in\Bb(X)$ by $K(x+y)\ :=-x+K_2y$ for $x\in V_2+V_3$ and $y\in X_2$. Then $K$ is compact and $(I+K)\lambda=\mu_2\subset\mu$. Thus $\lambda\sim^{gsc}\mu$ holds. 
\end{proof}

\begin{corollary}\label{c:semi-compact-reversible}
	Let $X$ be a Banach space
	with two closed linear subspaces $\lambda$ and $\mu$. Then the following hold.
	\begin{itemize}
		\item[(a)] We have 
		\begin{align}\label{e:semi-compact-reversible}
			(\lambda\sim^{sc}\mu,\mu\sim^{sc}\lambda)
			\Leftrightarrow&(\lambda\sim^c\mu)
			\Leftrightarrow (\lambda\sim^{sc}\mu, [\lambda-\mu]\in\ZZ),\\
			\label{e:global-semi-compact-reversible}
			(\lambda\sim^{gsc}\mu,\mu\sim^{gsc}\lambda)
			\Leftrightarrow&(\lambda\sim^{gc}\mu)
			\Leftrightarrow (\lambda\sim^{gsc}\mu, [\lambda-\mu]\in\ZZ).
		\end{align} 
		\item[(b)] Assume that $\lambda\overset{sc}\sim\mu$. Then we have
		\begin{align}\label{e:semi-compact-perturtation-ab-ba}
			[\lambda-\mu]=-[\mu-\lambda]. 
		\end{align}
	\end{itemize}
\end{corollary}

\begin{proof}
	1. By definition we have
	\begin{align*}(\lambda\sim^c\mu)
		&\Leftrightarrow (\lambda\sim^{sc}\mu, [\lambda-\mu]\in\ZZ),\\
		(\lambda\sim^{gc}\mu)
		&\Leftrightarrow (\lambda\sim^{gsc}\mu, [\lambda-\mu]\in\ZZ).
	\end{align*}
	\newline 2. If $\lambda\sim^{sc}\mu$ and $\mu\sim^{sc}\lambda$ hold, there are compact linear operators $K\in\Bb(\lambda,X)$ and $L\in\Bb(\mu,X)$ such that 
	\begin{align*}
		(I+K)\lambda\subset\mu\text{ and }(I+L)\mu\subset\lambda.
	\end{align*}
	By Lemma \ref{l:semi-compact-perturtation} we have 
	\begin{align*}
		0=[\lambda-\lambda]
		=[\lambda-\mu;K]+[\mu-\lambda;L].
	\end{align*}
	Then $I+K\colon\lambda\to\mu$ is Fredholm and $\lambda\sim^{c}\mu$ holds.
	\newline 3. If $\lambda\sim^{c}\mu$ holds, by Theorem \ref{t:semi-compact-perturbation} there are a closed linear subspace $\lambda_1$ of $\lambda$, a closed linear subspace $\mu_1$ of $\mu$ and a compact linear operator $K_1\in\Bb(\lambda,X)$ such that $I+K_1\colon\lambda_1\to\mu_1$ is a linear isomorphism.
	By Lemma \ref{l:semi-compact-perturtation} we have 
	\begin{align*}
		[\lambda-\mu]=[\lambda-\lambda_1]+[\lambda_1-\mu_1]+[\mu_1-\mu].
	\end{align*}
	Then we have $\dim\mu/\mu_1=-[\mu_1-\mu]\in\Z$. By Theorem \ref{t:semi-compact-perturbation} we have $\lambda\sim^{sc}\mu$ and $\mu\sim^{sc}\lambda$.
	\newline 4. Similar to Step 2 and Step 3, we have
	$(\lambda\sim^{gsc}\mu,\mu\sim^{gsc}\lambda)
	\Leftrightarrow(\lambda\sim^{gc}\mu)$. The only difference is that compact operators $K$ and $K_1$ there are defined on the whole space $X$.
	\newline 5. If $[\lambda-\mu]\in\ZZ$, \eqref{e:semi-compact-perturtation-ab-ba} follows from (a) and Lemma \ref{l:semi-compact-perturtation}.	
	If $[\lambda-\mu]\in\{+\infty,-\infty\}$, \eqref{e:semi-compact-perturtation-ab-ba} follows from Definition \ref{d:semi-compact-perturtation}.
\end{proof}

We have the following duality formula of relative dimension.

\begin{theorem}\label{t:relative-dimension-dual}
	Let $X$ be a Banach space with two closed linear subspaces $\lambda$ and $\mu$. Assume that $\lambda\sim^{gsc}\mu$. Then we have $\mu^\perp\sim^{gsc}\lambda^\perp$, and
	\begin{align}\label{e:relative-dimension-dual}
		[\mu^\perp-\lambda^\perp]=[\lambda-\mu].
	\end{align}
\end{theorem}

\begin{proof}
	1. Since $\lambda\sim^{gsc}\mu$, there is a compact linear operator $K\in\Bb(X)$ such that $(I_X+K)\lambda\subset\mu$. Then we have 
	$(I_{X^*}+K^*)\mu^\perp\subset\lambda^\perp$ and $\mu^\perp\sim^{gsc}\lambda^\perp$.
	\newline 2. By Theorem \ref{t:semi-compact-perturbation}, there is a compact linear operator $K_1\in\Bb(X)$, a closed linear subspace $\lambda_1$ of $\lambda$ and a closed linear subspace $\mu_1$ of $\mu$ such that $\dim\lambda/\lambda_1<+\infty$, $I_X+K_1\in\Bb(X)$ has bounded inverse and $(I_X+K_1)\lambda_1=\mu_1$. 
	By \cite[Theorem III.6.26]{Ka95}, $I_{X^*}+K_1^*\in\Bb(X)$ has bounded inverse.
	Since $(I_X+K_1)\lambda_1=\mu_1$, we have $(I_{X^*}+K_1^*)\mu_1^\perp=\lambda_1^\perp$. Then we have 
	\begin{align*}
		[\mu_1^\perp-\lambda_1^\perp]=0=[\lambda_1-\mu_1].
	\end{align*}
	\newline 3. By Lemma \ref{l:semi-compact-perturtation} we have
	\begin{align*}
		[\mu^\perp-\lambda^\perp]
		=&[\mu^\perp-\mu_1^\perp]+[\mu_1^\perp-\lambda_1^\perp]+[\lambda_1^\perp-\lambda^\perp]\\
		=&[\mu_1-\mu]+[\lambda_1-\mu_1]+[\lambda-\lambda_1]\\
		=&[\lambda-\mu].
	\end{align*}
\end{proof}

\subsection{An index formula for compact perturbations} \label{ss:stable-iindex-compact-perturbations} 
In this subsection we will prove an index formula for compact perturbations. 
For that, we
shall use the concepts of approximate nullity (approximate deficiency) defined by \cite[Section IV.4]{Ka95}:

\begin{definition}
	Let $M,N$ be closed linear subspaces of a Banach space $Z$.
	\newline (a) We define the \textit{approximate nullity} of the pair $M,N$,  denoted by $\nuli'(M,N)$, as the
	least upper bound of the set of integers $m$ ($m=+\infty$ being permitted) with the property that,
	for any $\e>0$, there is an $m$-dimensional closed linear subspace $M_{\e}\subset M$ with
	$\delta(M_{\e},N)<\e$.
	\newline (b) We define the \textit{approximate deficiency} of the pair $M,N$, denoted by $\defi'(M,N)$, by
	$\defi'(M,N):=\nuli'(M^{\bot},N^{\bot})$.
\end{definition}

\begin{note} While 
	\begin{align*}
		\nuli(M,N)\ :&=\dim M\cap N\text{ and }
		\defi(M,N)\ :=\dim Z/(M+N)
	\end{align*} 
	are defined in a
	purely algebraic fashion, the definitions of $\nuli^\prime(M,N)$ and $\defi^\prime(M,N)$ depend on the
	topology of the underlying space $Z$. Moreover, it is easy to show (see l.c., \cite[Theorems IV.4.18 and
	IV.4.19]{Ka95}) that
	\begin{align*}
		\nuli^\prime(M,N) \ &=\ \begin{cases} \nuli(M,N),&\text{for $M+N$ closed,}\\ +\infty,
			&\text{else,}\end{cases}%
		\tand\\%
		\defi^\prime(M,N)\ &=\ \begin{cases} \defi(M,N),&\text{for $M+N$ closed,}\\ +\infty, &\text{else.}\end{cases}
	\end{align*}
\end{note}

Then we have the following generalization of \cite[Theorem 5.26]{Ka95} and \cite[Proposition A.7.6.h]{BoZh18}.

\begin{theorem}\label{t:index-compact-perturbations}
	Let $X$ be a Banach space with two pairs $(M,N)$ and $(M^\prime,N^\prime)$ of closed linear subspaces. 
	\begin{itemize}
		\item [(a)] Assume that the pair $(M,N)$ is left Fredholm, $M^\prime\sim^{sc}M$ and $N^\prime\sim^{sc}N$. Then the pair $(M^\prime,N^\prime)$ is semi-Fredholm and we have
		\begin{align}\label{e:index-compact-perturbation}
			\Index(M^\prime,N^\prime)=\Index(M,N)
			+[M^\prime-M]+[N^\prime-N].
		\end{align}
		\item [(b)] Assume that the pair $(M,N)$ is right Fredholm, $M\sim^{gsc}M^\prime$ and $N\sim^{gsc}N^\prime$. Then the pair $(M^\prime,N^\prime)$ is right semi-Fredholm and \eqref{e:index-compact-perturbation} holds.
	\end{itemize}
\end{theorem}

\begin{proof} 
	The proof is similar to that of \cite[Theorem IV.5.26]{Ka95}.
	\newline (a) (i) By Theorem \ref{t:semi-compact-perturbation}, there are two compact linear operators $K_1\in\Bb(M_1^\prime,X)$ and $L_1\in\Bb(L_1^\prime,X)$, a closed linear subspace $M_1^\prime$ of $M^\prime$ and a closed linear subspace $N_1^\prime$ of $N^\prime$ such that
	\begin{align*}
		&\dim M^\prime/M_1^\prime<+\infty,\quad M_1\ :=(I+K_1)M_1^\prime\subset M,\quad\ker(I+sK_1)=\{0\},\\
		&\dim N^\prime/N_1^\prime<+\infty,\quad N_1\ :=(I+L_1)N_1^\prime\subset N,\quad\ker(I+sL_1)=\{0\}
	\end{align*} 
	for each $s\in[0,1]$. By \cite[Theorem IV.5.26]{Ka95}, $(I+sK_1)M_1^\prime$ and $(I+sL_1)N_1^\prime$ are closed for each $s\in[0,1]$, $M_1\ :=(I+K_1)M_1^\prime\subset M$ and $N_1\ :=(I+L_1)N_1^\prime\subset M$. By \cite[Lemma A.10]{WuZh25}, the families $\{(I+sK_1)M_1^\prime;\;s\in[0,1]\}$ and $\{(I+sL_1)N_1^\prime;\;s\in[0,1]\}$ are continuous. By Lemma \ref{l:subspace-index}, $(M_1,N_1)$ is left Fredholm and we have
	\begin{align}
		\Index(M_1,N_1)=\Index(M,N)-[M-M_1]-[N-N_1].
	\end{align} 
	
	(ii) Assume that $M\cap N=\{0\}$. Then we have 
	\begin{align}\label{e:compact-nuli-M1-N1}
		\nuli^\prime(M_1^\prime,N)<+\infty.
	\end{align} 
	
	To this end it is convenient to apply the result of \cite[Theorem IV.4.23]{Ka95}. Suppose that there is a sequence $u_n\in M_1^\prime$ such that $\|u_n\|=1$ and $\dist(u_n,N)\to 0$; we have to show that $\{u_n\}$ has a convergent subsequence. 
	Since $(I+K_1)^{-1}\in\Bb(M_1,M_1^\prime)$, $\{u_n\}$ is bounded.
    Since $K_1$ is compact, there is a subsequence $\{u_{n_j}\}$ of $\{u_n\}$ such that $K_1u_{n_j}\to w$. Since 
    \begin{align}\label{e:dist-compact-uMN}
    	\dist((I+K_1)u_n-K_1u_n,N)\to 0,
    \end{align}
    $(1+K_1)u_n\in M$ and $M+N$ is closed, we have $w\in\overline{M+N}=M+N$.
    Then $w=w_1+w_2$ for some $w_1\in M$ and $w_2\in N$. By \eqref{e:dist-compact-uMN} we have $\dist((I+K_1)u_{n_j}-w_1,N)\to 0$. Since $M+N$ is closed and $M\cap N=\{0\}$, by \cite[Theorem IV.4.2]{Ka95} we have $(I+K_1)u_{n_j}\to w_1$. Then we have $u_{n_j}\to(I+K_1)^{-1}w_1$. By \cite[Theorem IV.4.23]{Ka95}, we obtain \eqref{e:compact-nuli-M1-N1}.
    
    (iii) By (ii) and \cite[Theorems IV.4.18 and IV.4.19]{Ka95}, the pair $(M_1^\prime,N)$ is left semi-Fredholm if $M\cap N=\{0\}$. By \cite[Problem IV.4.6 and Problem IV.4.7]{Ka95}, the pair $(M^\prime,N)$ is left semi-Fredholm if $M\cap N=\{0\}$. By the proof of Lemma \ref{l:subspace-index}, the pair $(M^\prime,N)$ is left semi-Fredholm. 
    
    (iv) By (iii), the pair $(M^\prime,N^\prime)$ is left semi-Fredholm. By Lemma \ref{l:subspace-index}, the pair $(M_1,N_1)$ is left semi-Fredholm. Since $(I+sK_1)M_1^\prime\sim^c M_1$ and $(I+sL_1)N_1^\prime\sim^c N_1$, the pair $((I+sK_1)M_1^\prime,(I+sL_1)N_1^\prime)$is left semi-Fredholm. By Lemma and  \cite[Theorem IV.4.30 and Remark IV.4.31]{Ka95}, we have $\Index(M_1^\prime,N_1^\prime)=\Index(M_1,N_1)$.
    By Lemma \ref{l:subspace-index} and Lemma \ref{l:semi-compact-perturtation}, we have
    \begin{align*}
    	\Index(M^\prime,N^\prime)=&\Index(M_1^\prime,N_1^\prime)+[M^\prime-M_1^\prime]+[N^\prime-N_1^\prime]\\
    	=&\Index(M_1,N_1)+[M^\prime-M_1]+[N^\prime-N_1]\\
    	=&\Index(M,N)+[M^\prime-M]+[N^\prime-N].
    \end{align*}
    \newline (b) By (a) and Theorem \ref{t:relative-dimension-dual}.
\end{proof}

\subsection{Stability theorem for semi-compact perturbation of complemented linear subspaces }\label{ss:stable-semi-compact-perturbations-complemented}
The theory of semi-compact perturbation of complemented linear subspaces can be simplified. 

We have the following basic fact.

\begin{lemma}\label{l:subspace-operator} 
	
Let $X$ be a vector space.
Let $\alpha,\beta, \gamma$ be linear subspaces of $X$ such that $ \alpha\oplus\beta=\beta\oplus\gamma=X$ holds.
Then the following holds.
\newline (a) There is a linear operator $A\colon\alpha\rightarrow\beta$ such that
\begin{equation}\label{e:operator-space}
\gamma=\{a+Aa;\;x\in\alpha\}=\Graph(A),\quad\dom(A)=\alpha.
\end{equation}
Moreover, the linear operator $F=I_\alpha+A\colon\alpha\to\gamma$ is a linear isomorphism.
\newline (b) Assume that $X$ is a Banach space and $\alpha,\beta,\gamma$ are closed. Then the operator $A$ is bounded and the operator $F$ is bounded with bounded inverse.
\newline (c)
Assume that $X$ is a Banach space, $\alpha,\beta,\gamma$ are closed and there holds $\alpha\sim^{sc}\gamma$. Then the operator $A$ is compact, $\alpha\sim^{gc}\gamma$ and $[\alpha-\gamma]=0$.
\end{lemma}

\begin{proof} (a)
For each $c\in \gamma$, there are some $a\in\alpha,b\in \beta$ such that $c=a+b$ holds.
If there is $b^\prime\in\beta$ such that $a+b^\prime\in \gamma$, then $b^\prime-b\in \beta\cap\gamma=\{0\}$. Then $\gamma$ is the graph of a linear operator $A$.

Since $X=\beta\oplus\gamma$, for any $a\in \alpha$, $a=u+v$ with $u\in\beta,v\in\gamma$.
It follows that $v=a-u$.
So we have $\gamma=\Graph(A)$ with $\dom(A)=\alpha$. 

By the definition of $F$, the operator $F$ is a linear isomorphism.
\newline (b)
Since $X$ is a Banach space and $\alpha,\beta,\gamma$ are closed, by the closed graph theorem,
the operator $A$ is bounded and the operator $F$ is bounded with bounded inverse.
\newline (c)
Since $\alpha\sim^{sc}\gamma$, there is a compact linear operator $K\in \Bb(\alpha,X)$ such that $(I+K)\alpha\subset\gamma$.

Denote by $P$ the projection of $X$ along $\beta$ onto $\alpha$. For each $x\in\alpha$, we have $x+PKx\in\alpha$, $(I-P)Kx\in\beta$ and $x+PKx+(I-P)Kx=x+Kx\in\gamma$. Then we have 
\begin{align}\label{e:A-PK}
	A(I+PK)x=(I-P)Kx.
\end{align}
Since $PK\in\Bb(\alpha)$ is compact, $I_\alpha+PK\in\Bb(\alpha)$ is a Fredholm operator of index $0$. Then there are finite dimensional linear spaces $V_1$, $V_2$ of $\alpha$ and closed linear spaces $\alpha_1$, $\alpha_2$ of $\alpha$ such that
\begin{align*}
	\alpha&=V_1\oplus\alpha_1=V_2\oplus\alpha_2,\quad\dim V_1=\dim V_2,\\
	V_1&=\ker(I_\alpha+PK),\quad\alpha_2=\ran(I_\alpha+PK).
\end{align*}
By \eqref{e:A-PK} we have
\begin{align*}
	A|_{\alpha_2}=(I-P)K|_{\alpha_1}(I_{\alpha_1}+PK_{\alpha_1})^{-1}.
\end{align*}
Then the operator $A|_{\alpha_2}$ is compact. Since $\dim\alpha/\alpha_2<+\infty$ holds, the operator $A$ is compact. Since $\gamma=F\alpha$, by (a) we have $\alpha\sim^{gc}\gamma$ and $[\alpha-\gamma]=0$.
\end{proof}

We have the following nice result in the complemented case.

\begin{theorem}\label{t:stable-semi-compact-perturbations-complemented}
	Let $X$ be a Banach space. Let $(M,N)$ and be a pair of closed linear subspaces of $X$ with $M\sim^c N$. Assume that $M$ is complemented. Then there is a $\varepsilon>0$ such that for each pair $(M^\prime,N^\prime)$ of closed linear subspaces of $X$ with $M^\prime\sim^c N^\prime$, $\hat\delta(M^\prime,M)<\varepsilon$ and $\hat\delta(N^\prime,N)<\varepsilon$, we have $N\in\Ss^c(X)$,$M^\prime\in\Ss^c(X)$, $N^\prime\in\Ss^c(X)$, $M\sim^{gc} N$, $M^\prime\sim^{gc} N^\prime$, and
	\begin{align}\label{e:stable-semi-compact-perturbations-complemented}
		[M^\prime-N^\prime]=[M-N].
	\end{align}
\end{theorem}

\begin{proof}
	1. {\em We have $N\in\Ss^c(X)$, $M\sim^{gc} N$. Moreover, there are projections $P, Q\in\Bb(X)$ such that $\image P=M$ and $\image Q=N$ respectively, $P-Q$ is compact and $\ker P$ is a finite change of $\ker Q$.}
	
	Since $M\in\Ss^c(X)$, there is a linear subspace $\alpha_1\in\Ss(X)$ such that $X=M\oplus\alpha_1$.
	By Theorem \ref{t:index-compact-perturbations}, we have 
	\begin{align*}
		\Index(N,\alpha)=\Index(M,\alpha)+[N-M]=[N-M].
	\end{align*}
	Similar to the proof of \cite[Lemma A.2.6]{BoZh18}, there is a finite dimensional linear subspace $V$ of $X$, closed linear subspace $N_1$ of $N$ and closed linear subspace $\alpha_1$ of $\alpha$ such that 
	\begin{align*}
		N=N\cap \alpha\oplus N_1,\;\alpha=N\cap \alpha\oplus \alpha_1,\;X=V\oplus(N+\alpha).
	\end{align*}
	Then we have $N\cap\alpha_1=N\cap\alpha\cap\alpha_1=\{0\}$, and
	\begin{align*}
		X=N\cap\alpha\oplus N_1\oplus\alpha_1\oplus V.
	\end{align*}
	Then we have $N\in\Ss^c(X)$. 
	
	Since $X=M\oplus\alpha=M\oplus N\cap\alpha_1\oplus\alpha_1$, By Lemma \ref{l:semi-compact-perturtation}.b, we have $M\oplus N\cap\alpha\sim^c N_1\oplus V$ and $[(M\oplus N\cap\alpha)-(N_1\oplus V)]=0$. By Lemma \ref{l:subspace-operator}.c, we have $M\oplus N\cap\alpha\sim^{gc} N_1\oplus V$ and there is a compact linear operator $A:M\oplus N\cap\alpha\to\alpha_1$ such that $N_1\oplus V=\Graph(A)$. Then there are projections $P, Q\in\Bb(X)$ such that $\image P=M$ and $\image Q=N$ respectively, $P-Q$ is compact and $\ker P$ is a finite change of $\ker Q$.
	By \cite[Proposition A.7.6]{BoZh18}, $M\sim^{gc}N$ holds.
	\newline 2. By \cite[Lemma A.4.5]{BoZh18}, there is a $\varepsilon>0$ such that for each pair $(M^\prime,N^\prime)$ of closed linear subspaces of $X$ with $M^\prime\sim^c N^\prime$, $\hat\delta(M^\prime,M)<\varepsilon$ and $\hat\delta(N^\prime,N)<\varepsilon$, we have 
	\begin{align*}
		X=M^\prime\oplus\ker P=N^\prime\oplus\ker Q.
	\end{align*} 
	Then we have $M^\prime\in\Ss^c(X)$ and $N^\prime\in\Ss^c(X)$. By Step 1, we have $M^\prime\sim^{gc}N^{\prime}$. Denote by $P^\prime$ the projection on $M^\prime$ along $\ker P$ and $Q^\prime$ the projection on $N^\prime$ along $\ker Q$ respectively. Since $\ker P$ is a finite change of $\ker Q$, by Lemma \ref{l:subspace-operator}.c, the operator $Q^\prime-P^\prime$ is compact. 
	Since $\varepsilon$ is small enough, we have $\|P^\prime-P\|<1$ and $\|Q^\prime-Q\|<1$. Similar to the proof of \cite[Lemma 2.4]{ZhLo99}, by \cite[Lemma 1.4.10]{Ka95}, there are linear maps $S,T\in\Bb(X)$ with bounded inverse such that
	\begin{align*}
		SP=P^\prime S,\quad TQ=Q^\prime T
	\end{align*}
	Moreover, we have $\lim\limits_{P^\prime\to P}S=I$ and $\lim\limits_{Q^\prime\to Q}T=I$.
	By \cite[Theorem IV.5.26]{Ka95}, we have
	\begin{align*}
		[M^\prime-N^\prime]=&[\Index Q^\prime P^\prime\colon\image P^\prime\to\image Q^\prime]\\
		=&[\Index TQT^{-1} SPS^{-1}\colon\image SP\to\image TQ]\\
		=&[\Index QT^{-1} SP\colon\image P\to\image Q]\\
		=&[\Index QP\colon\image P\to\image Q]=[M-N].
	\end{align*}
\end{proof}

\subsection{An isometric linear embedding of a Banach space into some Banach space with (AP)} As indicated in the introduction, the proof of Theorem \ref{t:stable-relative-dimention} is rather complicated due to the the lack of control of the norm $\|K^\prime-K\|$, where $K$ and $K^\prime$ are compact operators such that $(I+K)M\subset N$ and $(I+K^\prime)M^\prime\subset N^\prime$. In order to overcome the difficulty, we firstly embed the Banach space $X$ into some Banach space with (AP). 

A. Grothendieck \cite{Groth55} defined the following. 

\begin{definition}\label{d:Banach-AP} 
	(a) A Banach space $X$ has the {\em approximation property} if the following condition is satisfied:
	\begin{itemize}
		\item[(AP)] For every compact subset 
		$K$ and every $\varepsilon> 0$, there exists a finite rank operator $T$ 
		such that $\|x-Tx\|\le c$ whenever $x\in K$.
	\end{itemize}  
	
	(b) A Banach space $X$ has the {\em bounded approximation property} if we can find
	a constant $c\ge 1$ such that the following condition is satisfied:
	\begin{itemize}
		\item[(BAP)] For every compact subset 
		$K$ and every $\varepsilon> 0$, there exists a finite rank operator $T$ 
		of norm not greater than $c$ such that $\|x-Tx\|\le c$ whenever $x\in K$.
	\end{itemize}  
	
	(c) A Banach space $X$ has the {\em metric approximation property (MAP)} if the constant $c$ can be chosen to be $1$ in (b).
\end{definition} 

Recall the notion of approximation property for compact operators.

\begin{definition}\label{d:ap}
	Let $X$, $Y$ be two Banach spaces. 
	We call a compact operator $K\in\Bb(X,Y)$ has {\em approximation property (AP)}, if for each $\varepsilon>0$, there is a finite rank operator $F\in\Bb(X,Y)$ such that $\|K-F\|<\varepsilon$. 	
\end{definition}

A. Grothendieck \cite[Chapter I, p. 165]{Groth55} discovered many properties that are
equivalent to (AP), one of them is 
\begin{itemize}
	\item[(AP1)] For every Banach space $Y$, a compact operator $T\in\Bb(Y,X)$ can be  approximated by finite rank operators. 
\end{itemize} 

P. Enflo \cite{Enflo73} showed that $l_p$ with $2<p\le+\infty$ has a subspace without (AP). Later, A. Szankowski  subsequently showed that $l_p$ with $1\le p<2$ has a subspace without (AP) in \cite{Sz78} and $\Bb(H)$ has not (AP) for a infinite dimensional Hilbert space $H$ \cite{Sz79}.

The Banach-Mazur theorem says that every separable Banach space over $\KK$ is isometrically isomorphic to a closed subspace of $C([0,1],\KK)$. For the general case, we have the following result.

\begin{lemma}\label{l:embed-into-AP}	
	Let $X$ be a Banach space. Denote by $\Kk$ the unit ball of $X^*$, equipped with the $w^*$-topology. Define $j\colon X\to C(\Kk)$ by 
	\begin{align}\label{e:definition-j}
		j(x)(x^*)=x^*(x).
	\end{align}
	Then $\Kk$ is a compact Hausdorff space, $j$ is an isometric linear embedding, and $C(\Kk)$ has (MAP).  
\end{lemma}

\begin{proof} 
	By Banach-Alaoglu theorem, $\Kk$ is a compact Hausdorff space. The map $j$ is linear, and it is isometric by Hahn-Banach theorem. By A. Grothendieck \cite[Chapter I, p. 185]{Groth55}, $C(\Kk)$ has (MAP).
\end{proof}

\subsection{Stability theorem for global semi-compact perturbation }\label{ss:stable-semi-compact-perturbations}
As indicated in the introduction, the proof of Theorem \ref{t:stable-relative-dimention} is rather complicated due to the the lack of control of the norm $\|K^\prime-K\|$, where $K$ and $K^\prime$ are compact operators such that $(I+K)M\subset N$ and $(I+K^\prime)M^\prime\subset N^\prime$. 

We need some preparation to prove Theorem \ref{t:stable-relative-dimention}. 

By \cite[Lemma 2.24]{LhcZ24}, we have

\begin{lemma}\label{l:operator-same-space}
	Let $X$, $Y$ be two Banach spaces. Let $A\in\Bb(X,Y)$ be a bounded operator with bounded $C=A^{-1}|_{AX}:AX\to X$, $B\in\Bb(X,Y)$ be a bounded operator. 
	Then the following hold. 
	\begin{itemize}
		\item[(a)] The space $AX$ is closed, and we have
		\begin{align}\label{e:operator-same-space}
			\delta(AX,BX)\le \|C\|\|A-B\|.
		\end{align}
		\item[(b)] Assume that there hold $A^{-1}\in\Bb(AX,X)$ and $\|A-B\|\|A^{-1}\|<1$. Then we have $B^{-1}\in\Bb(BX,X)$, $\|B^{-1}\|\le(\|A^{-1}\|^{-1}-\|A-B\|)^{-1}$ and $BX$ is closed.		
	\end{itemize}	
\end{lemma}

\begin{proof}[Proof of Theorem \ref{t:stable-relative-dimention}]	
	(a) follows from Lemma \ref{l:semi-compact-perturtation}.
	\newline (b) By Theorem \ref{t:semi-compact-perturbation}, there are closed subspaces $M_1$ of $M$ and $N_1$ of $N$ and a compact linear operator $K_1\in\Bb(X)$ such that 
	\begin{align}\label{e:sc-finite-codimension-app1}
		\dim M/M_1<+\infty,
		\text{ and }
		(I+K_1)M_1=N_1.
	\end{align}        
	Moreover, the operators $I+K_1$ and $I+K_1\colon M_1\to N_1$ are linear isomorphisms. 
	Since \eqref{e:sc-finite-codimension-app1}, $N_1\subset(I+K_1)M\cap N$ and 
	$[M-N]\in\Z$ hold, 
	there are finite dimensional linear subspaces $V$ of $(I+K_1)M$ and $U$ of $N$ respectively such that
	\begin{align*}
		V\oplus(I+K_1)M\cap N&=(I+K_1)M,\\			
		U\oplus(I+K_1)M\cap N&\begin{cases}
			=N &\text{ if }[M-N]\in\ZZ,\\
			\subset N&\text{ if }[M-N]=-\infty,
		\end{cases}\\
		\dim U&=\dim V+m \text{ if }[M-N]=-\infty.
	\end{align*}     
	Then we obtain
	\begin{align}\label{e:finite-change-same-space-proof}
		(I+K_1)M\oplus U\begin{cases}
			=N\oplus V &\text{ if }[M-N]\in\ZZ,\\
			\subset N\oplus V&\text{ if }[M-N]=-\infty.
		\end{cases}
	\end{align}
    \newline (c)    
    (i) By \cite[(IV.4.6) and Theorem IV.4.2]{Ka95}, we have $\gamma(V,N)\ge\frac{\gamma(N,V)}{1+\gamma(N,V)}>0$. By \cite[(1.4.1)]{Neu68} we have
    \begin{align}\label{e:delta-N+V}
    	\delta(N+V,N^\prime+V)\le\delta(N,N^\prime)\gamma(N,V)^{-1}.
    \end{align}
    
    Since $I+K_1$ has bounded inverse, by Lemma \ref{l:operator-same-space} we have
    \begin{align}
    	\label{e:com-iso-gap}
    	\delta((I+K_1)M^\prime,(I+K_1)M)\le a\delta(M^\prime,M).
    \end{align}
    Since there hold $U\cap (I+K_1)M=\{0\}$ and
    \begin{align*}
    	\delta(M^\prime, M)<a^{-1}\gamma(U,(I+K_1)M),
    \end{align*}    
    by \cite[Problem IV.4.15, Theorem IV.4.19 and Theorem IV.4.24]{Ka95} we have 
    \begin{align*}
    	(I+K_1)M^\prime\cap U=\{0\}.
    \end{align*}    
    By \cite[(IV.4.34)]{Ka95} and \eqref{e:com-iso-gap}, we have 
    \begin{align*}
    	\gamma(U,(I+K_1)M^\prime)\ge\frac{\gamma(U,(I+K_1)M)-a\delta(M^\prime,M)}{1+a\delta(M^\prime,M)}>0.
    \end{align*}   
    By \cite[(1.4.1)]{Neu68}, we have 
    \begin{align}\label{e:delta-M-prime+U}
    	\delta((I+K_1)M^\prime+U,(I+K_1)M+U)
    	\le\theta_1.  	
    \end{align}    
    By Lemma \ref{l:delta-MNL}, \eqref{e:finite-change-same-space-proof}, \eqref{e:delta-N+V} and \eqref{e:delta-M-prime+U} we have
    \begin{align}\label{e:sum-gap}
    	\delta((I+K_1)M^\prime+U,N^\prime+V)\le \eta_1-1<1.
    \end{align}	
    \newline (ii)
    By Lemma \ref{l:semi-compact-perturtation} and \cite[Proposition A.7.6.a and Proposition A.7.6.e]{BoZh18}, we have
    \begin{align*}
    	(I+K_1)M^\prime+U\overset{sc}\sim N^\prime+V.
    \end{align*}
    By Theorem \ref{t:semi-compact-perturbation}, there is a compact operator $K_2\in\Bb((I+K_1)M^\prime\oplus U,X)$ such that either
    \begin{align*}
    	(I+K_2)((I+K_1)M^\prime\oplus U)\subset N^\prime+V
    \end{align*}
    or a compact operator $K_2\in\Bb(N^\prime+V,X)$
    \begin{align*}
    	(I+K_2)(N^\prime+V)\subset (I+K_1)M^\prime+ U.
    \end{align*}
    
    Denote by $j\colon X\to C(\Kk)$ the isometric linear embedding obtained in Lemma \ref{l:embed-into-AP}.    
    Since $C(\Kk)$ has (AP), there is a finite rank linear operator $F\in\Bb(\dom(K_2), C(\Kk))$ with $\|jK_2-F\|<\frac{2-\eta_1}{2}$. Set $K^\prime_2\ :=jK_2-F$. By the construction of $K^\prime_2$ and Lemma \ref{l:finite-change}, we have
    \begin{align}\label{e:finite-difference-W1-NV}
    	W_1\sim^f j(N^\prime+V)\text{ or }W_2\sim^f j((I+K_1)M^\prime+ U),
    \end{align}  
     where
    \begin{align*}
    	W_1\ :=&(j+K^\prime_2)((I+K_1)M^\prime+ U),\\
    	W_2\ :=&(j+K^\prime_2)(N^\prime+V).
    \end{align*}  

    Since $K^\prime_2-jK_2$ is an operator of finite rank, the operator $K^\prime_2$ is compact. Since $\|K^\prime_2\|<\frac{2-\eta_1}{2}<1$, by Lemma \ref{l:semi-compact-perturtation}.d.
    we have
    \begin{align}\label{e:relative-dimension-W1-MU}
    	[j((I+K_1)M^\prime+U)-W_1]=[W_2-j(N^\prime+V)]=0.
    \end{align}
    Since $\|K^\prime_2\|<\frac{2-\eta_1}{2}$, by Lemma \ref{l:operator-same-space}, we have
    \begin{align}\label{e:delta-W1-M+U}
    	\delta(W_1,j((I+K_1)M^\prime+U))<\frac{2-\eta_1}{\eta_1},\quad \delta(j(N^\prime+V),W_2)<\frac{2-\eta_1}{\eta_1}.
    \end{align}
    By Lemma \ref{l:delta-MNL}, \eqref{e:sum-gap} and \eqref{e:delta-W1-M+U},
    we have
    \begin{align}\label{e:delta-small-W1-NV}
    	\delta(W_1,j(N^\prime+V))<1,\quad \delta(j((I+K_1)M^\prime+U),W_2)<1.
    \end{align}
    \newline (iii)    
    By \cite[Lemma A.3.1.d]{BoZh18} and \eqref{e:delta-small-W1-NV}, we have
    \begin{align*}
    	\delta(\frac{W_1}{W_1\cap j(N^\prime+V)},\frac{j(N^\prime+V)}{W_1\cap j(N^\prime+V)})<&1,\text{ or}\\
    	\delta(\frac{j((I+K_1)M^\prime+U)}{W_2\cap j((I+K_1)M^\prime+U)},\frac{W_2}{W_2\cap j((I+K_1)M^\prime+U)})<&1.
    \end{align*}
    By \cite[Corollary IV.2.6]{Ka95} and \eqref{e:finite-difference-W1-NV}, we have
    \begin{align*}
    	\dim\frac{W_1}{W_1\cap j(N^\prime+V)}\le&\dim\frac{j(N^\prime+V)}{W_1\cap j(N^\prime+V)},\text{ or}\\
    	\dim\frac{j((I+K_1)M^\prime+U)}{W_2\cap j((I+K_1)M^\prime+U)}\le&\dim\frac{W_2}{W_2\cap j((I+K_1)M^\prime+U)}.
    \end{align*}    
    Then we obtain
	\begin{align}\label{e:relative-dimension-W1-NV}
		[W_1-j(N^\prime+V)]\le 0,\quad\text{ or}\quad [j((I+K_1)M^\prime+U)-W_2]\le 0.
	\end{align}
	\newline (iv)
	By Step (i), Lemma \ref{l:relative-dimension-embedding}, Lemma \ref{l:semi-compact-perturtation}, \eqref{e:relative-dimension-W1-MU} and \eqref{e:relative-dimension-W1-NV} we have
	\begin{align*}
		[M^\prime-N^\prime]=&
		[jM^\prime-jN^\prime]\\
		=&[jM^\prime-j((I+K_1)M^\prime+U)]+[j((I+K_1)M^\prime+U)-W_i]\\
		&+[W_i-j(N^\prime+V)]+[j(N^\prime+V)-jN^\prime]\\
		\le&-\dim U+\dim V\\
		=&\begin{cases}
			[M-N] &\text{if }[M-N]\in\ZZ,\\
			-m &\text{if }[M-N]=-\infty
		\end{cases}
	\end{align*}
    for $i=1$ or $i=2$.
	\newline(d) If $[M-N]=-\infty$, we are done.
	
	If $[M-N]\in\ZZ$, by Corollary \ref{c:semi-compact-reversible} we have $N\sim^{sc} M$ and $[M-N]=-[N-M]$. By (c) we have
	\begin{align*}
		[M^\prime-N^\prime]=-[N^\prime-M^\prime]\ge-[N-M]=[M-N].
	\end{align*}
    \newline(e) By (c) and (d).
\end{proof}

\begin{lemma}\label{l:meta}
	Let $J$ be a connected topological space and $f\colon J\to \Z\cup\{+\infty,-\infty\}$ be a continuous function. Then $f$ is a constant function.
\end{lemma}

\begin{proof}
	For each $n\in\Z$ we have
	\begin{align*}
		&f^{-1}((-\infty,n]\cup\{-\infty\})=J\setminus f^{-1}((n,+\infty)\cup\{+\infty\}),\\
		&f^{-1}([n,+\infty)\cup\{+\infty\})=J\setminus
		f^{-1}((-\infty,n)\cup\{-\infty\}),\\
		&f^{-1}(n)=f^{-1}((-\infty,n]\cup\{-\infty\})\cap f^{-1}([n,+\infty)\cup\{+\infty\}).
	\end{align*}
	Then $f^{-1}((-\infty,n]\cup\{-\infty\})$, $f^{-1}([n,+\infty)\cup\{+\infty\})$ and $f^{-1}(n)$ are three open and closed subsets of $J$. Since $J$ is connected, each of the three subsets is either $J$ or an emptyset. 
	
	If $f^{-1}(n)=J$ we are done. Assume that $f^{-1}(n)\ne J$ for all $n\in\Z$. Then we have $f^{-1}(\Z)=\emptyset$. Consequently, $f^{-1}(-\infty)=f^{-1}((-\infty,n]\cup\{-\infty\})$ and $f^{-1}(+\infty)=f^{-1}([n,+\infty)\cup\{+\infty\})$ are two open and closed subsets of $J$, and $J=f^{-1}(-\infty)\cup f^{-1}(+\infty)$ is a disjoint union. Since $J$ is connected, there holds either $J=f^{-1}(-\infty)$ or $J=f^{-1}(+\infty)$.
\end{proof}

\begin{proof}[Proof of Theorem \ref{t:stable-relative-dimention-family}]
	By Theorem \ref{t:stable-relative-dimention}, the function $f(s)$ defined by 
	\begin{align*}
		f(s)=[M(s)-N(s)],\quad\forall s\in J
	\end{align*} 
	is a continuous function. By Lemma \ref{l:meta}, the function $f$ is a constant function on $s\in J$. 	 	 
\end{proof}

\section{Stability theorems for semi-Fredholm tuples}\label{s:stable-Fredholm}

We need some preparations for the proof of Theorem \ref{t:stable-semi-Fredholm}, Theorem \ref{t:stable-Fredholm-index} and Theorem \ref{t:stable-Fredholm-index-family}.

Firstly we have the following version of triangular inequality of gap.
 
\begin{lemma}\label{l:delta-MNL}
	Let $X$ be a Banach space with three linear subspaces $M$, $N$, $L$. Then we have
	\begin{align}\label{e:delta-MNL}
		1+\delta(M,L)\le(1+\delta(M,N))(1+\delta(N,L)).
	\end{align}
\end{lemma}

\begin{proof} 
	We note that the gap can be defined for any two linear subspace of $X$.  
	By \cite[Lemma IV.2.2]{Ka95} we have
	\begin{align}\label{e:delta-u-NL}
		\begin{split}
			\|u\|+\dist(u,L)&\le(\|u\|+\dist(u,N))(1+\delta(N,L))\\
			&\le\|u\|(1+\delta(M,N))(1+\delta(N,L)).
		\end{split}		
	\end{align}
    for all $u\in M$. By the definition of the gap, we obtain \eqref{e:delta-MNL}.
\end{proof}

By the similar argument we have

\begin{lemma}\label{l:delta-MNL-uv}
	Let $X$ be a Banach space with three closed linear subspaces $M$, $N$, $L$. Assume that $M\supsetneq L$. Then for each pair of positive numbers $(a, b)$ with $(1+a)(1+b)<2$, either there exists a $u\in \romS_M$ with $\dist(u,N)>a$ or there exists a $v\in \romS_N$ with $\dist(v,L)>b$. 
\end{lemma}

\begin{proof} 
	By \cite[Lemma III.1.2]{Ka95}, there exists a $u\in\romS_M$ with
	\begin{align*}
		\dist(u,L)>(1+a)(1+b)-1.
	\end{align*}   
	By \cite[Lemma IV.2.2]{Ka95} we have
	\begin{align*}		
		(1+a)(1+b)<1+\dist(u,L)\le(1+\dist(u,N))(1+\delta(N,L)).
	\end{align*}
	Assume that $\dist(u,N)\le a$. Then we have $\delta(N,L)>b$. By the definition of the gap, there exists a $v\in \romS_N$ with $\dist(v,L)>b$. 
\end{proof}

We have the following generalization of \cite[Lemma A.3.2]{BoZh18}.

\begin{lemma}\label{l:ball-distances}
    Let $X$ be a Banach space over $\KK$ with a closed linear subspace $M$. Let $u_1,\ldots,u_n\in \operatorname{S}_X$ be on the unit sphere of $X$. Set%
	\[
	V_k\ :=\begin{cases} \{0\},& \text{for $k=0$},\\
		\Span\{u_1,\dots,u_k\}, & \text{for $k=1\ldots,n$}.
	\end{cases}
	\]
	Assume that $\dist(u_k,M+V_{k-1})\ge \delta_k$ for $k=1,\ldots,n-1$ and positive numbers $\delta_k$. Set 
	\begin{align*}
		\Delta_k\ :=\frac{\prod_{i=k}^n\delta_i}{\prod_{i=k+1}^n(1+\delta_i)}
	\end{align*}
	for $k=1\ldots,n$. Then we
	have $\delta_k\le 1$, $\dim (M+V_k)/M=k$, and
	\begin{align}\label{e:ball-distances}
		 &\|\dist(\sum_{k=1}^{n}a_ku_k,M)\| \ge \frac{\Delta_1}{n}\sum_{k=1}^n|a_k|,\\
		 \label{e:ball-distances-gamma}
		 &\gamma(V_n,M)\ge \frac{\Delta_1}{n}.
	\end{align}
	for all $a_1,\ldots,a_k\in\KK$. 
\end{lemma}

\begin{proof} 	
	Since $\delta_k>0$, we have $u_k\notin
	M+V_{k-1}$. Then we have $\dim (M+V_k)/(M+V_{k-1})=1$, and by induction we have $\dim (M+V_k)/M=k$.
	
	We have $\dist(a_1u_1,M)\ge\delta_1|a_1|$, and 
	\begin{align*}
		\dist(\sum_{i=1}^{k}a_iu_i,M)&\ge  \max\{\delta_k|a_k|,\dist(\sum_{i=1}^{k-1}a_iu_i,M)-|a_k|\}\\
		&\ge\ \max\biggl\{\frac{\delta_k}{1+\delta_k}\dist(\sum_{i=1}^{k-1}a_iu_i,M),\delta_k|a_k|\biggr\}
	\end{align*}
    for $k=2,\ldots,n$. By induction we have
	\begin{align*}
		\dist(\sum_{k=1}^{n}a_ku_k,M)&\ge
		\max\biggl\{\Delta_k|a_k|;\;
		 k=1,\ldots,n\biggr\}\\
		&\ge\ \frac{\Delta_1}{n}\sum_{k=1}^n|a_k|.
	\end{align*}
    Since $\|\sum_{k=1}^{n}a_ku_k\|\le\sum_{k=1}^n|a_k|$ holds, we have $\gamma(V_n,M)\ge\frac{\Delta_1}{n}$.
\end{proof}

For the gap between two finite dimensional linear subspace, we have the following lemma.

\begin{lemma}\label{l:gap-finite-dimension}
	Let $X$ be a Banach space over $\KK$ with two closed linear subspaces $M$, $M^\prime$. Let  $\{\delta_k;\;k=1,\ldots,n \}\subset(0,1]$ be a set of positive numbers. Let
	$v_1,\ldots,v_n\in \romS_X$ be such that%
	\begin{align*}
		\dist(v_k,M+V_{k-1}) \ge\delta_k 
	\end{align*}
	for
	\begin{align*}
		V_k\ :=
		\begin{cases} 
			\{0\},& \text{for $k=0$},\\
			\Span\{v_1,\dots,v_k\},  &\text{for $k=1\ldots,n$}.
		\end{cases}
	\end{align*}
	Let $\delta\ge 0$ be a positive number. 
	Set
	\begin{align*}
		\Delta\ :=\frac{\prod_{k=1}^n\delta_k}{n\prod_{k=2}^n(1+\delta_k)}.
	\end{align*}
    Let $v_k^\prime\in X$ for each $k=1,\ldots,n$ be such that 
	\begin{align}\label{e:condition-gap-finite-dimension}
		\|v_k-v_k^\prime\|\le\delta<\frac{\Delta-\delta(M^\prime,M)}{1+\delta(M^\prime,M)}.
	\end{align}
    Set
    \begin{align*}
    	V_k^\prime\ :=
    	\begin{cases} 
    		\{0\},& \text{for $k=0$},\\
    		\Span\{v_1^\prime,\dots,v_k^\prime\},  &\text{for $k=1\ldots,n$}.
    	\end{cases}
    \end{align*}
	Then we have $\dim V_n=\dim V_n^\prime=n$, and
	\begin{align}
		\label{e:gamma-finite-extension}
		\gamma(V_n^\prime,M^\prime)&\ge\frac{\Delta-\delta(M^\prime,M)-(1+\delta(M^\prime,M))\delta}{(1+\delta)(1+\delta(M^\prime,M))},\\
		\label{e:gap-finite-dimension}
		\hat\delta(V_n,V_n^\prime)&\le \frac{\delta}{\Delta-\delta}.
	\end{align}
\end{lemma}

\begin{proof}
	Let $(a_1,\ldots,a_n)\in\KK^n$ be arbitrary. By Lemma \ref{l:ball-distances}, we have
	\begin{align*}
		\Delta\sum_{k=1}^n|a_k|\le\dist(\sum_{k=1}^{n}a_kv_k,M)\le\|\sum_{k=1}^{n}a_kv_k\|\le\sum_{k=1}^n|a_k|.
	\end{align*}
	
	Then we have $\dim (M+V_n)/M=n$. Note that we have
	\begin{align*}
		\|\sum_{k=1}^n a_kv_k^\prime\|\le (1+\delta)\sum_{k=1}^n|a_k|.
	\end{align*}
	By \cite[Lemma IV.2.2]{Ka95} we have
	\begin{align}\label{e:lower-bound-norm-v-prime}
		\begin{split}
			\dist(\sum_{k=1}^n a_kv_k^\prime,M^\prime)
			\ge&\frac{\dist(\sum_{k=1}^{n}a_kv_k^\prime,M)-\sum_{k=1}^n\|a_kv_k^\prime\|\delta(M^\prime,M)}{1+\delta(M^\prime,M)}\\
			\ge& \frac{\Delta-\delta-(1+\delta)\delta(M^\prime,M)}{1+\delta(M^\prime,M)}\sum_{k=1}^n|a_k|.
		\end{split}				
	\end{align}
	By \eqref {e:lower-bound-norm-v-prime}, we have $\dim (M^\prime+V_n^\prime)/M^\prime=n$ and \eqref{e:gamma-finite-extension}.
	Apply \eqref{e:lower-bound-norm-v-prime} for $M=M^\prime=\{0\}$, we obtain
	\begin{align*}
		\sum_{k=1}^n \|a_kv_k^\prime\|
		\ge (\Delta-\delta)\sum_{k=1}^n|a_k|.
	\end{align*}
    Note that
	\begin{align*}
		\sum_{k=1}^n\|a_k(v_k-v_k^\prime)\|\le
		\delta\sum_{k=1}^n|a_k|.
	\end{align*}
    By the definition of the gap we obtain \eqref{e:gap-finite-dimension}.
\end{proof}

We have the following existence of nearby finite dimensional linear subspace in a nearby closed linear subspace (cf. \cite[Lemma 1.7]{Neu68}). 

\begin{lemma}\label{l:existence-finite-dimension-embeding-close}
	Let $X$ be a Banach space with two closed linear subspaces $N$ and $N^\prime$. Let $V$ be a linear subspace of $N$ with $\dim V=n<+\infty$. Assume that $\delta(N,N^\prime)<\frac{1}{2^{n-1}n}$. Then for each $\varepsilon\in(0,\frac{1}{2^{n-1}n}-\delta(N,N^\prime))$, there is a linear subspace $V^\prime$ of $N^\prime$ such that there hold $\dim V^\prime=n$ and
	\begin{align}\label{e:existence-finite-dimension-embeding-close}
		\hat\delta(V,V^\prime)\le\frac{2^{n-1}n(\delta(N,N^\prime)+\varepsilon)}{1-2^{n-1}n(\delta(N,N^\prime)+\varepsilon)}.
	\end{align}
\end{lemma}

\begin{proof} 
	By induction and \cite[Lemma IV.2.3]{Ka95}, there exist
	$v_1,\ldots,v_n\in \romS_N$ such that%
	\begin{align*}
	    \dist(v_k,V_{k-1}) =1 
	\end{align*}
	for
	\begin{align*}
	    V_k\ :=
	    \begin{cases} 
	    	\{0\},& \text{for $k=0$},\\
		    \Span\{v_1,\dots,v_k\},  &\text{for $k=1\ldots,n$}.
	    \end{cases}
    \end{align*}
	Then $V_n=V$. Moreover, there exists a $v_k^\prime\in N^\prime$ for each $k=1,\ldots,n$ with 
	\begin{align*}
		\|v_k-v_k^\prime\|<\delta_\varepsilon\ :=\delta(N,N^\prime)+\varepsilon.
	\end{align*}
    Set $V^\prime\ :=V_n^\prime$. By Lemma \ref{l:gap-finite-dimension}, we are done.
\end{proof}

\begin{proposition}\label{p:stable-finite-extension}
	Let $X$ be a Banach space. Let $(M,N)$ be a pair of closed linear subspaces with $M\subset N$. Let $V$ be a closed linear subspace $V$ of $N$ with $M\cap V=\{0\}$ such that the space $M+V$ is closed. Then the following hold.
	\begin{itemize}
		\item [(a)] Assume that $N=M\oplus V$. Then for each pair $(M^\prime,N^\prime)$ of closed linear subspaces with $M^\prime\subset N^\prime$ and
		\begin{align}\label{e:stable-finite-extension-1-condition}			
			(1+\delta(N^\prime,N))(1+\delta(M,M^\prime)\gamma(M,V)^{-1})<2,
		\end{align}  
	    we have
		\begin{align}\label{e:stable-finite-extension-1}
			\dim N^\prime/M^\prime\le \dim N/M.
		\end{align}
	    \item [(b)] For each pair $(M^\prime,N^\prime)$ of closed linear subspaces with $M^\prime\subset N^\prime$ and
	    \begin{align}\label{e:stable-finite-extension-2-condition}\begin{split}
	    		&\delta(M^\prime, M)< \gamma(V,M),\\
	    		&(1+\frac{(1+\gamma(V,M))\delta(M^\prime,M)}{\gamma(V,M)-\delta(M^\prime,M)})(1+\delta(N,N^\prime))<2,
	    	\end{split}	    	
	    \end{align}
        We have
        \begin{align}\label{e:stable-finite-extension-2}
        	\dim N^\prime/M^\prime\ge\dim V.
        \end{align} 
    \item [(c)] 
        We have
        \begin{align}\label{e:stable-finite-extension}
        	\dim N^\prime/M^\prime=\dim N/M.
        \end{align} 
        if both the conditions in (a) and (b) hold.
	\end{itemize}
\end{proposition}

\begin{proof}
	By \cite[(IV.4.6) and Theorem IV.4.2]{Ka95}, we have $\gamma(V,M)\ge\frac{\gamma(M,V)}{1+\gamma(M,V)}>0$. By \cite[(1.4.1)]{Neu68} we have
	\begin{align}\label{e:delta-M+V}
		\delta(M+V,M^\prime+V)\le\delta(M,M^\prime)\gamma(M,V)^{-1}.
	\end{align}
    
    Assume that $\delta(M^\prime, M)<\gamma(V,M)$. Since $V\cap M=\{0\}$ holds, by \cite[Problem IV.4.15, Theorem IV.4.19 and Theorem IV.4.24]{Ka95} we have    
    We have $M^\prime\cap V=\{0\}$ and $M^\prime+V$ is closed.
    By \cite[(IV.4.34)]{Ka95}, we have 
    \begin{align}\label{e:gamma-V-M-prime}
    	\gamma(V,M^\prime)\ge\frac{\gamma(V,M)-\delta(M^\prime,M)}{1+\delta(M^\prime,M)}>0.
    \end{align}   
    By \cite[(1.4.1)]{Neu68}, we have 
    \begin{align}\label{e:delta-M-prime+V}
    	\delta(M^\prime+V,M+V)\le\frac{(1+\gamma(V,M))\delta(M^\prime,M)}{\gamma(V,M)-\delta(M^\prime,M)}.  	
    \end{align}    
    \newline (a) The inequality \eqref{e:stable-finite-extension-1} holds if $\dim N/M=+\infty$. Assume that we have $\dim N/M<+\infty$. Then the linear subspace $V$ with $M\oplus V=N$ is finite dimensional and closed. By Lemma \ref{l:delta-MNL} and \eqref{e:delta-M+V} we have
    \begin{align*}
    	1+\delta(N^\prime,M^\prime+V)&\le (1+\delta(N^\prime,N))(1+\delta(N,M^\prime+V))\\
    	&\le(1+\delta(N^\prime,N))(1+\delta(M,M^\prime)\gamma(M,V)^{-1}).
    \end{align*}
    By \cite[Lemma A.3.1.d]{BoZh18} we have \begin{align*}
    	\delta(N^\prime/M^\prime,(M^\prime+V)/M^\prime)=\delta(N^\prime,M^\prime+V)<1.
    \end{align*}
    By \cite[Corollary IV.2.6]{Ka95} we have
    \begin{align*}
    	\dim N^\prime/M^\prime\le\dim V=\dim N/M.
    \end{align*}
    \newline (b) Since $M\cap V=\{0\}$ and $\delta(M^\prime, M)< \gamma(V,M)$, we have $M^\prime\cap V=\{0\}$. Since $\delta(M+V,N)=0$, by Lemma \ref{l:delta-MNL} and \eqref{e:delta-M-prime+V} we have
    \begin{align*}
    	1+&\delta(M^\prime+V,N^\prime)\le (1+\delta(M^\prime+V,M+V))(1+\delta(N,N^\prime))\\
    	&\le(1+\frac{(1+\gamma(V,M))\delta(M^\prime,M)}{\gamma(V,M)-\delta(M^\prime,M)})(1+\delta(N,N^\prime)).
    \end{align*}
    By \cite[Lemma A.3.1.d]{BoZh18} we have \begin{align*}
    	\delta((M^\prime+V)/M^\prime,N^\prime/M^\prime)=\delta(M^\prime+V,N^\prime)<1.
    \end{align*}
    By \cite[Corollary IV.2.6]{Ka95} we have
    \begin{align*}
    	\dim N^\prime/M^\prime\ge\dim V.
    \end{align*}
    \newline (c) By (a) and (b).
\end{proof}

\begin{theorem}[Splitting theorem]\label{t:existence-finite-dimension-embeding-close-oposite}
	Let $X$ be a Banach space with four closed linear subspaces $M$, $N$, $L$ and $S$ with $M=L\oplus S$ and $\dim S=n<+\infty$. Let $a\in(0,\sqrt{2}-1)$ be a positive number. For a nonnegative integer $i$, we set
	\begin{align*}
		a_i\ :&=\begin{cases}
			\frac{a^i}{i(1+a)^{i-1}},&\text{if }i\ge 1,\\
			1&\text{if } i=0,
		\end{cases}\\
	    	\delta_i\ :&=\begin{cases}
	    	\min\{\frac{\delta(N,M)(1+a_i)}{a_i},1\},&\text{if }i\ge 1,\\
	    	\delta(N,M)&\text{if } i=0.
	    \end{cases}
	\end{align*} 
    \begin{itemize}
    	\item[(a)] Assume that there holds
    	\begin{align}\label{e:conditions-existence-V}
    		(1+\delta_{n+1})(1+\delta(L,N)\gamma(L,S)^{-1})<2.
    	\end{align} 
    	Then there is a linear subspace $V_k$ of $M$ 
    	such that 
    	\begin{align}\label{e:properties-existence-V-dimension}   		
    		&\dim V_k=k\le n,\quad N\cap V_k=\{0\},\\
    		\label{e:properties-existence-V-MN}   			
    		&\gamma(V_k,N)\ge a_k, \quad\delta(M,N+V_k)\le a,\quad\delta(N+V_k,M)\le\delta_k.    	 	
    	\end{align} 
        \item[(b)] Assume that there hold  \eqref{e:conditions-existence-V} and 
        \begin{align}\label{e:conditions-existence-V-more}
        	\delta(L,N)<a_k.
        \end{align} 
        Then $V_k$ can be chosen such that
        \begin{align}\label{e:properties-existence-V-more-LV} 
        	&L\cap V_k=\{0\},\\
        	\label{e:properties-existence-V-more-gamma}           	
        	&\gamma(V_k,L)\ge\frac{ a_k-\delta(L,N)}{1+\delta(L,N)},\\	
        	\label{e:properties-existence-V-more-delta} 
        	&\delta(L+V_k,N+V_k)\le\frac{(1+a_k)\delta(L,N)}{a_k-\delta(L,N)}.        	  		
        \end{align}
        \item[(c)] Assume that there hold  \eqref{e:conditions-existence-V}, \eqref{e:conditions-existence-V-more} and
        \begin{align}\label{e:conditions-existence-V-c}
        	c_k\ :=\frac{(n-k)(a+1)^{n-k-1}(1+\delta_k)\delta_k}{(a-\delta_k)^{n-k}(1-\delta_k)}\in[0,1).
        \end{align}  
        Then there are two linear subspaces $U_{n-k}$, $V_k$ of $M$ and a linear subspace $W_{n-k}$ of $N$ respectively such that (a), (b) hold and  
        \begin{align}
        	\label{e:properties-existence-V-more-LVW} 
        	&(L+V_k)\cap W_{n-k}=\{0\},\\
        	\label{e:splitting-M}
        	&M=L\oplus V_k\oplus U_{n-k},  \quad\dim U_{n-k}=n-k,\\
        	\label{e:gamma-U-LV}
        	&\gamma(U_{n-k},L+V_k) 
        	\begin{cases}
        		>\frac{(a-\delta_k)^{n-k}}{(n-k)(a+1)^{n-k-1}(1+\delta_k)},&\text{if }k<n,\\
        		=1&\text{if }k=n,
        	\end{cases}\\
        	\label{e:delta-U-W}
        	&\hat\delta(U_{n-k},W_{n-k})
        	\le\frac{c_k}{1-c_k}.
        \end{align}
    \end{itemize} 	
\end{theorem}

\begin{proof} 
	1. Let $v_1,\ldots, v_{k_1}\in \romS_M$ be on the unit sphere of $M$ with $\dist(v_{i},N+V_{i-1})> a$, where 
	\begin{align*}
		V_i\ :=
		\begin{cases} 
			\{0\},& \text{for $i=0$},\\
			\Span\{v_1,\dots,v_i\},  &\text{for $i=1\ldots,k_1$}.
		\end{cases}
	\end{align*} 
    Then we have $\dim(N+V_{k_1})/N=\dim V_{k_1}=k_1\le n$.
    
	In fact, for each $i=1,\ldots, k_1$, by Lemma \ref{l:ball-distances}, we have $\dim(N+V_i)/N=i$, and 
	\begin{align*}
		\gamma(V_i,N)\ge a_i.
	\end{align*} 
    By \cite[(1.4.1)]{Neu68} we have
	\begin{align*}		
		\delta(N+V_i,M)=&\delta(V_i+N,M+M)\\
		\le&\min\{\frac{\delta(N,M)(1+\gamma(V_i,N))}{\gamma(V_i,N)},1\}\le\delta_i.
	\end{align*}  
    Assume that $k_1>n$. Then we can take $i=n+1$. By \eqref{e:conditions-existence-V} we have 
	\begin{align*}		
		(1+\delta(N+V_{n+1},M))(1+\delta(L,N)\gamma(L,S)^{-1})<2.
	\end{align*}
    By Proposition \ref{p:stable-finite-extension}.a we have 
    \begin{align*}
    	n+1=\dim(N+V_{n+1})/N\le \dim M/L=n,
    \end{align*} 
    a contradiction.
    
    We take $k$ to be the maximal integer of all such $k_1$. Then we have $\delta(M,N+V_k)\le a$. Thus $V_k$ satisfies our requirements and (a) is proved.
	\newline 2. By Step 1 we have 
	\begin{align*}
		N\cap V_k=\{0\},\quad\delta(L,N)<a_k\le\gamma(V_k,N).
	\end{align*}
    By \cite[Theorem IV.4.24]{Ka95} we have $L\cap V_k=\{0\}$. By \cite[Lemma 1.4]{Neu68} we have
    \begin{align*}
    	\gamma(V_k,L)\ge&\frac{ \gamma(V_k,N)-\delta(L,N)}{1+\delta(L,N)}\\	 
    	\ge&\frac{ a_k-\delta(L,N)}{1+\delta(L,N)},\\	
    	\delta(L+V_k,N+V_k)\le&\delta(L,N)(1+\gamma(V_k,L)^{-1})\\
    	\le&\frac{(1+a_k)\delta(L,N)}{a_k-\delta(L,N)}.
    \end{align*}
    Then (b) is proved.
	\newline 3.	Let $\varepsilon\in(0,1-\delta(N+V_k,M))$ be a positive number. Set $\delta_\varepsilon\ :=\delta(N+V_k,M)+\varepsilon$. 
	We proceed by induction with nonnegative integers $j$ with $j\le n-k$ to construct linear subspaces $U_j$ of $M$ and $W_j$ of $N+V_k$. 
	
	Set $W_0\ :=\{0\}$. 
	Assume that $j<n-k$ holds and linear spaces $U_j\subset M$, $W_j\subset N$ is given respectively with
	\begin{align*}
		U_j&=\Span\{u_1,\ldots,u_j\},\\		
		W_j&=\Span\{w_1,\ldots,w_j\},
	\end{align*}
	where $u_1,\ldots,u_j\in M$, $w_1,\ldots,w_j\in \romS_{N+V_k}$, and 
	\begin{align*}
		\|u_s-w_s\|<\delta_\varepsilon,&\quad s=1,\ldots,j,\\
		\dist(w_s,L+V_k+U_{s-1})>b\ :=\sqrt{2}-1,&\quad s=0,\ldots,j
	\end{align*} 	
	By Step 1 we have $\delta(M,N+V_k)\le a$. By Lemma \ref{l:delta-MNL-uv}, there exists a $w_{j+1}\in \romS_{N+V_k}$ with $\dist(w_{j+1},L+V_k+U_j)>b$.
	Then there exists a $u_j\in M$ with
	\begin{align*}
		\|u_{j+1}-w_{j+1}\|<\delta_\varepsilon.
	\end{align*}  
    Set
    \begin{align*}
    	U_{j+1}&\ :=\Span\{u_1,\ldots,u_{j+1}\},\\		
    	W_{j+1}&\ :=\Span\{w_1,\ldots,w_{j+1}\}.
    \end{align*}
    Then the required linear spaces $U_j$ and $W_j$ for $j=0,\ldots, n-k$ is defined.
    \newline 4. We now prove (c). Choose sufficiently small $\varepsilon>0$. Set $\delta_\varepsilon\ :=\delta(N+V_k,M)+\varepsilon$. Then we have $\delta_\varepsilon\le a+\varepsilon<b$. Note that we have
    \begin{align*}
    	&\dist(u_j,L+V_k)\ge\dist(u_j,L+V_k+U_{j-1})>b-\delta_\varepsilon>0,\\
    	&0<1-\delta_\varepsilon<\|u_j\|<1+\delta_\varepsilon
    \end{align*}	
    for $j=1,\ldots,n-k$. Then we obtain
    \begin{align}
    	\nonumber
    	\dist(u_j,U_{j-1})&\ge\dist(u_j,L+V_k+U_{j-1})\\
    	\label{e:dist-uLVU}
    	&>\frac{b-\delta_\varepsilon}{1+\delta_\varepsilon}\|u_j\|\ge\frac{b-\delta_\varepsilon}{1+\delta_\varepsilon}\dist(u_j,L+V_k),\\
    	\label{e:dist-uw}
    	\|u_j-w_j\|&<\frac{\delta_\varepsilon}{1-\delta_\varepsilon}\|u_j\|,\\
    	\label{e:dist-uw-LV}
    	\dist(u_j-w_j,L+V_k)&<\frac{\delta_\varepsilon(1+\delta_\varepsilon)}{b-\delta_\varepsilon}\dist(u_j,L+V_k)
    \end{align}	
    for $j=1,\ldots,n-k$.
    By Lemma \ref{l:ball-distances}, \eqref{e:dist-uLVU} and Step 2 we have $M=L\oplus V_k\oplus U_{n-k}$, $\dim U_{n-k}=n-k$, and
    \begin{align*}
    	\gamma(U_{n-k},L+V_k)&\ge\frac{(\frac{b-\delta_\varepsilon}{1+\delta_\varepsilon})^{n-k}}{(n-k)(1+\frac{b-\delta_\varepsilon}{1+\delta_\varepsilon})^{n-k-1}}\\ 
    	&>\frac{(\frac{a-\delta_k}{1+\delta_k})^{n-k}}{(n-k)(1+\frac{a-\delta_k}{1+\delta_k})^{n-k-1}}\\ 
    	&=\frac{(a-\delta_k)^{n-k}}{(n-k)(a+1)^{n-k-1}(1+\delta_k)}.
    \end{align*}
    By Lemma \ref{l:gap-finite-dimension}, \eqref{e:dist-uLVU}, \eqref{e:dist-uw} and \eqref{e:dist-uw-LV}, we have $\dim W_{n-k}=n-k$, $(L+V_k)\cap W_{n-k}=\{0\}$, and
    \begin{align*}
    	\hat\delta(U_{n-k},W_{n-k})
    	&\le\frac{\frac{\delta_\varepsilon}{1-\delta_\varepsilon}}{\frac{(\frac{b-\delta_\varepsilon}{1+\delta_\varepsilon})^{n-k}}{(n-k)(1+\frac{b-\delta_\varepsilon}{1+\delta_\varepsilon})^{n-k-1}}-\frac{\delta_\varepsilon}{1-\delta_\varepsilon}}\\
    	&\le\frac{\frac{\delta_k}{1-\delta_k}}{\frac{(\frac{a-\delta_k}{1+\delta_k})^{n-k}}{(n-k)(1+\frac{a-\delta_k}{1+\delta_k})^{n-k-1}}-\frac{\delta_k}{1-\delta_k}}\\   	
    	&=\frac{c_k}{1-c_k}.
    \end{align*}
\end{proof}

Here we define the index of a semi-Fredholm tuple.

\begin{definition}\label{d:index-semi-Fredholm-tuple}
	Let $X$ be a vector space with four linear subspaces $Y_1$, $M$, $N$ and $Y_2$ such that $Y_1\subset M\cap N\subset M+N\subset Y_2$. 
	\begin{itemize}
		\item [(a)] We define 
		\begin{align}\label{e:index-semi-Fredholm-techad}
			\Index(Y_1;M,N;Y_2)\ :=\dim(M\cap N/Y_1)-\dim Y_2/(M+N) 
		\end{align}
		if either $\dim(M\cap N/Y_1)$ or $\dim Y_2/(M+N)$ is finite.
		\item [(b)] The tuple $(Y_1;M,N;Y_2)$ is called {\em Fredholm} if both $\dim(M\cap N/Y_1)$ and $\dim Y_2/(M+N)$ are finite. 
	\end{itemize}
\end{definition}

\begin{definition}\label{d:semi-Fredholm-tuplle}
	Let $X$ be a Banach space with four closed linear subspaces $Y_1$, $M$, $N$ and $Y_2$ such that $Y_1\subset M\cap N\subset M+N\subset Y_2$. 
	The tuple $(Y_1;M,N;Y_2)$ is called  {\em semi-Fredholm} if either $\dim(M\cap N/Y_1)$ or $\dim Y_2/(M+N)$ are finite, and $M+N$ is closed. 	    
\end{definition}

\begin{proof}[proof of Theorem \ref{t:stable-semi-Fredholm}]
	We divide the proof into five steps.
	\newline 1. {\em Equation \eqref{e:stable-cap} in (b) holds.}
	
	By \cite[Lemma IV.4.4]{Ka95} and let $u\in M^\prime\cap N^\prime$ be arbitrary, we have
	\begin{align}\label{e:delta-cap-MN}
		\delta(M^\prime\cap N^\prime,M\cap N)\le\frac{2}{\gamma(M,N)}(\delta(M^\prime,M)+\delta(N^\prime,N)).
	\end{align}
	By Proposition \ref{p:stable-finite-extension}.a,  \eqref{e:stable-cap} holds if
	$Y_1^\prime\subset M^\prime\cap N^\prime$ with sufficiently small $\delta(Y_1,Y_1^\prime)$, $\delta(M^\prime,M)$ and
	$\delta(N^\prime,N)$. 
	\newline 2. {\em Assume that $M\cap N=Y_1$. Then $M^\prime+N^\prime$ in (b) is closed.}

    Since $Y_1^\prime\subset M^\prime\cap N^\prime$, by Step 1 we have $Y_1^\prime=M^\prime\cap N^\prime$. By \cite[Lemma 1.4]{Neu68}, $M^\prime+N^\prime$ is closed if
    $Y_1^\prime=M^\prime\cap N^\prime$ with sufficiently small $\delta(Y_1,Y_1^\prime)$, $\delta(M^\prime,M)$ and
    $\delta(N^\prime,N)$.
    \newline 3. {\em Assume that $M+ N=Y_2$. Then $M^\prime+N^\prime=Y_2^\prime$ holds in (a).}
    
    Apply Step 1 and Step 2 by consider the tuple $(Y_2^\perp;M^\perp,N^\perp)$ with the perturbation $((Y_2^\prime)^\perp;(M^\prime)^\perp,(N^\prime)^\perp)$. By \cite[Theorem IV.2.9 and Theorem IV.4.8]{Ka95}, the result follows.
    \newline 4. {\em The statement (a) holds.}
    
    There is a linear subspace $S$ of $Y_2$ with $Y_2=(M+N)\oplus S$ and $\dim S=n<+\infty$.
    
    Let $a\in(0,\sqrt{2}-1)$ be a positive number. For a nonnegative integer $i$, we set
    \begin{align*}
    	a_i\ :&=\begin{cases}
    		\frac{a^i}{i(1+a)^{i-1}},&\text{if }i\ge 1,\\
    		1&\text{if } i=0,
    	\end{cases}\\
    	\delta_i\ :&=\begin{cases}
    		\min\{\frac{\delta(Y_2^\prime,Y_2)(1+a_i)}{a_i},1\},&\text{if }i\ge 1,\\
    		\delta(Y_2^\prime,Y_2)&\text{if } i=0,
    	\end{cases}\\
        c_i\ :&=\frac{(n-i)(a+1)^{n-i-1}(1+\delta_i)\delta_ii}{(a-\delta_i)^{n-i}(1-\delta_i)}.
    \end{align*} 
    Since $\delta(M,M^{\prime})$, $\delta(N,N^{\prime})$ and $\delta(Y_2^{\prime},Y_2)$
    is sufficiently small, by Theorem \ref{t:existence-finite-dimension-embeding-close-oposite}, there are two linear subspaces $U$, $V$ of $Y_2$ and a linear subspace $W$ of $Y_2^\prime$ respectively such that there hold
    \begin{align*}
    	&Y_2=(M+N)\oplus V\oplus U, \quad\dim U=\dim W,\\
    	&k\ :=\dim V\le n,\quad Y_2^\prime\cap V=(M+N+V)\cap W=\{0\},\\		
    	&\gamma(V,Y_2^\prime)\ge a_k,\quad\delta(Y_2,Y_2^\prime+V)\le a, \quad\delta(Y_2^\prime+V,Y_2)\le\delta_k,\\  &\gamma(V,M+N)\ge\frac{ a_k-\delta(M+N,Y_2^\prime)}{1+\delta(M+N,Y_2^\prime)},\\	       	   	
    	&\gamma(U,M+N+V) 
    	\begin{cases}
    		>\frac{(a-\delta_k)^{n-k}}{(n-k)(a+1)^{n-k-1}(1+\delta_k)},&\text{if }k<n,\\
    		=1&\text{if }k=n,
    	\end{cases}\\
    	&0\le c_k<1,\quad\hat\delta(U,W)
    	\le\frac{c_k}{1-c_k}.
    \end{align*}
    By \cite[(1.4.1)]{Neu68} we have
    \begin{align*}
    	\delta(M+N,Y_2^\prime)&\le\delta(N,Y_2^\prime)+(\delta(M,Y_2^\prime)+\delta(N,Y_2^\prime))\gamma(M,N)^{-1}\\
    	&\le\delta(N,N^\prime)+(\delta(M,N^\prime)+\delta(N,N^\prime))\gamma(M,N)^{-1}.
    \end{align*}
    Since $V\cap(M+N)=\{0\}$, by \cite[(1.4.1)]{Neu68} we have
    \begin{align*}
    	\delta(M+V,M^\prime+V)&\le\delta(M,M^\prime)(
    	1+\gamma(V,M)^{-1})\\
    	&\le\delta(M,M^\prime)(
    	1+\gamma(V,M+N)^{-1}).
    \end{align*} 
    Since $U\cap(M+N+V)=\{0\}$, by \cite[(1.4.1)]{Neu68} we have
    \begin{align*}
        &\delta(M+V+U,M^\prime+V+W)\\
        \le&\delta(M+V,M^\prime+V)+(\delta(M+V,M^\prime+V)\\&
        +\delta(U,W))\gamma(U,M+V)^{-1}\\
        \le&\delta(M+V,M^\prime+V)+(\delta(M+V,M^\prime+V)\\
        &+\hat\delta(U,W))\gamma(U,M+N+V)^{-1}.
    \end{align*}    
    Since $(M+V+U)+N=Y_2$, by Step 3 we have $(M^\prime+V+W)+N^\prime=Y_2^\prime+V$. Then we have
    \begin{align*}
    	\dim Y_2^\prime/(M^\prime+N^\prime)&=\dim (Y_2^\prime+V)/(M^\prime+N^\prime)-\dim V\\
    	&=\dim(M^\prime+N^\prime+V+W)/(M^\prime+N^\prime)-\dim V\\
    	&\le\dim(V+W)-\dim V\\
    	&=\dim W\le\dim Y_2/(M+N).
    \end{align*}
    By \cite[Problem IV.4.7]{Ka95}, the subspace $M^\prime+N^\prime$ is closed.
    \newline 5. {\em The statement (b) holds.}
    
    Apply Step 4 by consider the tuple $(M^\perp,N^\perp;Y_1^\perp)$ with the perturbation $((M^\prime)^\perp,(N^\prime)^\perp;(Y_1^\prime)^\perp)$. By \cite[Theorem IV.2.9 and Theorem IV.4.8]{Ka95}, the result follows.    
\end{proof}

\begin{proof}[proof of Theorem \ref{t:stable-Fredholm-index}]
	(a) We start to continue from the end of Step 4 in the proof of Theorem \ref{t:stable-semi-Fredholm}. Then the subspace $M^\prime+N^\prime$ is closed. By Proposition \ref{p:finite-extension} and Theorem \ref{t:stable-semi-Fredholm}.a we have
	\begin{align*}
		&\Index((Y_1^\prime;M^\prime,N^\prime; Y_2^\prime)\\
		=&\Index((Y_1^\prime;M^\prime+V+W,N^\prime; Y_2^\prime+V)-\dim W\\
		=&\dim((M^\prime+V+W)\cap N^\prime)/Y_1^\prime-\dim W\\
		\ge&\dim((M+V+U)\cap N)/Y_1-\dim W\\
		=&\Index((Y_1;M+V+U,N; Y_2)-\dim U\\
		=&\Index((Y_1;M,N; Y_2)+\dim V\\ 
		\ge&\Index((Y_1;M,N; Y_2). 
	\end{align*}
	\newline (b) Apply (a) by consider the tuple $(Y_2^\perp;M^\perp,N^\perp;Y_1^\perp)$ with the perturbation $((Y_2^\prime)^\perp;(M^\prime)^\perp,(N^\prime)^\perp;(Y_1^\prime)^\perp)$. By \cite[Theorem IV.2.9 and Theorem IV.4.8]{Ka95}, the result follows.
	\newline (c) follows from (a) and (b). 
	\newline
	(d) We start to continue from the end of Step 4 in the proof of Theorem \ref{t:stable-semi-Fredholm}. Then the subspace $M^\prime+N^\prime$ is closed.
	Note that $\dim((M+V+W)\cap N)/Y_1=+\infty$.  By Proposition \ref{p:finite-extension} and Theorem \ref{t:stable-semi-Fredholm}.a we have
	\begin{align*}
		&\Index((Y_1^\prime;M^\prime,N^\prime; Y_2^\prime)\\
		=&\Index((Y_1^\prime;M^\prime+V+W,N^\prime; Y_2^\prime+V)-\dim W\\
		=&\dim((M^\prime+V+W)\cap N^\prime)/Y_1^\prime-\dim W\\
		\ge&m+\dim W-\dim W=m. 
	\end{align*}
	\newline (e) Apply (d) by consider the tuple $(Y_2^\perp;M^\perp,N^\perp;Y_1^\perp)$ with the perturbation $((Y_2^\prime)^\perp;(M^\prime)^\perp,(N^\prime)^\perp;(Y_1^\prime)^\perp)$. By \cite[Theorem IV.2.9 and Theorem IV.4.8]{Ka95}, the result follows.	
\end{proof}

\begin{proof}[proof of Theorem \ref{t:stable-Fredholm-index-family}] 
	By Theorem \ref{t:stable-Fredholm-index}, the function $f(s)$ defined by 
	\begin{align*}
		f(s)=\Index(Y_1(s);M(s),N(s);Y_2(s)),\quad\forall s\in J
	\end{align*} 
    is a continuous function. By Lemma \ref{l:meta}, the function $f$ is a constant function on $s\in J$. 	 	 
\end{proof}

\section{Perturbed augmented Morse index}\label{ss:stable-augmented-morse}

In this section we study the augmented Morse index of perturbed symmetric forms. For the corresponding study of the Morse index, see \cite[Section 2.3]{LhcZ24}. 

Let $V$ be a complex space over $\KK$ and $Q\colon V\times V\to\KK$ be a symmetric form.
For each subset $\lambda$ of $V$, we denote by 
\begin{align*}
	\lambda^Q\ :=\{u\in V;\;Q(u,v)=0,\;\forall v\in\lambda\}.
\end{align*}
Then $m^{\pm}(Q)$ and $m^0(Q)$ denote the {\em Morse positive (or negative) index} and the {\em nullity} of $Q$ respectively.  
For each {\em isotropic subspace} $\epsilon$ of $V$, i.e. $\epsilon$ is a linear subspace of $V$ with $\epsilon\subset \epsilon^Q$, we have the {\em reduced form} $\tilde Q$ on $\epsilon^Q/\epsilon$ defined by 
\begin{align}\label{e:reduced-symmetric-form}
	\tilde Q(x+\epsilon,y+\epsilon)\ :=Q(x,y),\quad\text{for all }x,y\in\epsilon^Q.
\end{align}  
We denote by $\pi_\epsilon\colon \epsilon^Q\to \epsilon^Q/\epsilon$ the natural projection.

We have the following two classical results on symmetric forms.

\begin{lemma}\label{l:Q-decomposition}
	Let $V$ be a vector space over $\KK$ and $Q\colon V\times V\to\KK$ be a symmetric form. Let $\alpha$ be a finite dimensional positive (negative) definite subspace of $V$ with respect to $Q$. Then we have the following direct sum decomposition
	\begin{align}\label{e:Q-decomposition}.
		V=\alpha\oplus\alpha^Q.
	\end{align}	
\end{lemma}

\begin{proof}
	Since $\alpha$ is a positive (or negative) definite subspace of $V$ with respect to $Q$, we have $\alpha\cap\alpha^Q=0$. Take an orthonormal bases $e_1,\ldots,e_k$ of $\alpha$ with respect to $\pm Q$. Then for each $v\in V$, we have 
	\begin{align*}
		v=\pm\sum_{j=1}^{k}Q(v,e_j)e_j+\left(v\mp\sum_{j=1}^{k}Q(v,e_j)e_j\right)\in\alpha+\alpha^Q.
	\end{align*} 
	Thus \eqref{e:Q-decomposition} holds. 
	
\end{proof}

\begin{lemma}\label{l:maximal-Q}
	Let $V$ be a vector space over $\KK$ and $Q\colon V\times V\to\KK$ be a symmetric form. Let $\alpha$ be a maximal positive (or negative) definite subspace of $V$ with respect to $Q$. Then 
	the form $Q|_{\alpha^Q}$ is negative (or positive) semi-definite.
\end{lemma}

\begin{proof}
	Since $\alpha$ is a positive (or negative) definite subspace of $V$ with respect to $Q$, we have $\alpha\cap\alpha^Q=0$. If the form $Q|_{\alpha^Q}$ is not negative (or positive) semi-definite, there is a $v\in \alpha^Q\setminus\alpha$ such that $\pm Q(v,v)>0$. Then $Q$ is positive (or negative) definite on the space $Q\oplus\Span\{v\}$. This contradicts to the fact that $\alpha$ is a maximal positive (or negative) definite subspace of $V$ with respect to $Q$. Thus the form $Q|_{\alpha^Q}$ is negative (or positive) semi-definite.
\end{proof}

We recall the following quantities of  symmetric forms defined by \cite[Definition 2.14]{LhcZ24}.

\begin{definition}\label{d:bounded-symmetric-pair} \cite[Definition 2.14]{LhcZ24}
	(a) Let $X$ be a Banach space over $\KK$ with a closed linear subspace $V$. Let $Q\colon V\times V\to\KK$ be a symmetric form. If 
	\begin{align}\label{e:norm-Q}
		\|Q\|\ :=\begin{cases}
			\sup\limits_{x,y\in V\setminus\{0\}}\frac{|Q(x,y)|}{\|x\|\|y\|}<+\infty&\text{ if } V\ne\{0\},\\0&\text{ if } V=\{0\},
		\end{cases}
	\end{align}
	we call $(Q,V)$ a {\em bounded symmetric pair} of $X$. 
	\newline (b) Let $X$ be a Banach space with a closed linear subspace $V$. Let $Q\colon V\times V\to\KK$ be a bounded semi-definite symmetric form. The {\em reduced minimum modulus} $\gamma(Q)$ of $Q$ is defined by 
	\begin{align}\label{e:reduced-minimum-modulus-Q}
		\gamma(Q)\ :=\begin{cases}
			\inf\limits_{x\in V\setminus V^Q}\frac{|Q(x,x)|}{(\dist(x,V^Q))^2}&\text{ if } V\ne V^Q,
			\\0&\text{ if } V=V^Q.
		\end{cases}
	\end{align}
	\newline (c) Let $X$ be a Banach space with two symmetric pairs $(Q,V)$ and $(R,W)$. Let $c\ge 0$ be a real number. We define the {\em $c$-gap} $\delta_c(Q,R)$ between $R$ and $Q$ to be the infimum of the non-negative number $\delta$ such that 
	\begin{align}\label{e:c-gap-pair}
		\begin{split}
			|Q(x,y)&-R(u,v)|\le \delta(\|u\|+\|x\|)(\|v\|+\|y\|)+\\
			& c\left((\|u\|+\|x\|)\|v-y\|+\|u-x\|(\|v\|+\|y\|)\right)
		\end{split}    	
	\end{align}
	for all $x,y\in V$ and $u,v\in W$.
\end{definition}

\begin{lemma}\label{l:V-Q-closed} 
	Let $X$ be a Banach space with a bounded symmetric pair $(Q,V)$. Then $V^Q$ is closed.
\end{lemma}

\begin{proof}
	Since $V$ is closed and $Q$ is bounded, for each $x\in V$, $x^Q\ :=\{y\in V;\;Q(x,y)=0\}$ is closed. Then $V^Q=\bigcap\limits_{x\in V}x^Q$ is closed.
\end{proof}

The following lemma estimate the gap between annihilators of finite dimensional definite linear subspaces.

\begin{lemma}\label{l:annihilator-close}
	Let $X$ be a Banach space with two bounded symmetric pairs $(Q,V)$ and $(R,W)$. Let $\alpha\subset V$ and $\beta\subset W$ be two closed linear subspaces of $X$. Let $h\in\{1,-1\}$ and $c\ge 0$ be two real numbers. Set
	\begin{align}\label{e:eta}
		\eta\ :=\gamma(Q|_\alpha)^{-1}(2+\delta(\alpha,\beta))((2+\delta(\alpha,\beta))\delta_c(Q,R)+ 2c\delta(\alpha,\beta)).
	\end{align}
    Assume that $\dim\alpha=\dim\beta=n<+\infty$ holds, $hQ$ is positive definite on $\alpha$ and $\eta<n^{-1}$. Then $hR$ is positive definite on $\beta$, and we have
    \begin{align}\label{e:annihilator-close}
    	\delta(\alpha^Q,\beta^R)\le&\delta(V,W)
    	+n(1-n\eta)^{-\frac{1}{2}}(1+\delta(\alpha,\beta))\gamma(Q|_\alpha)^{-\frac{1}{2}}C,\\ 
    	\label{e:annihilator-close-gamma}
    	\gamma(R|_\beta)\ge& n^{-1}(1-n\eta)(1+\delta(\alpha,\beta))^{-1}\gamma(Q|_\alpha).
    \end{align}
    where 
    \begin{align*}
    	C\ :=&\gamma(Q|_\alpha)^{-\frac{1}{2}}\delta_c(Q,R)(2+\delta(V,W))((1-n\eta)^{-\frac{1}{2}}(1+\delta(\alpha,\beta))+1)\\
    	&+
    	c\gamma(Q|_\alpha)^{-\frac{1}{2}}\left((2+2\delta(V,W))(1-n\eta)^{-\frac{1}{2}}(1+\delta(\alpha,\beta))-2\right).
    \end{align*} 
\end{lemma}

\begin{proof}
	1. Take an orthonormal basis $\{v_1,\ldots,v_n\}$ of $(\alpha,hQ|_\alpha)$ and an $\varepsilon>0$. 
    By \eqref{e:reduced-minimum-modulus-Q} we have
    \begin{align*}
    	\|v_i\|\le\gamma(Q|_\alpha)^{-\frac{1}{2}}.
    \end{align*}
    Then there is a subset $\{u_i\}_{i=1,\ldots,n}$ of $\beta$ such that
    \begin{align*}
    	\|v_i-u_i\|\le(1+\varepsilon)\delta(\alpha,\beta)\|v_i\|.
    \end{align*}
    holds for each $i=1,\ldots,n$. Then we have
    \begin{align}\label{e:norm-u}
    	\max\{\|u_i\|;\;i=1,\ldots,n\}\le(1+(1+\varepsilon)\delta(\alpha,\beta))\gamma(Q|_\alpha)^{-\frac{1}{2}}.
    \end{align}
    \newline 2. Set $A\: =\left(hR(u_i,u_j)\right)_{i,j=1,\ldots,n}$.
    For each matrix $M=\left(a_{ij}\right)_{i,j=1,\ldots,n}$, we set
    \begin{align}\label{e:norm-sum-colum}
    	\|M\|_1\ :=\max_j\{\sum_{i=1}^{n}|a_{ij}|\}.
    \end{align}
    By \eqref{e:c-gap-pair} we have 
    \begin{align*}
    	|(A-I_n)_{i,j}|\le\eta_\varepsilon
    	\ :=&\delta_c(Q,R)\gamma(Q|_\alpha)^{-1}(2+(1+\varepsilon)\delta(\alpha,\beta))^2\\ &+c(1+\varepsilon)\delta(\alpha,\beta)\gamma(Q|_\alpha)^{-1}(4+2(1+\varepsilon)\delta(\alpha,\beta))
    \end{align*}
    for $i,j=1,\ldots,n$.
    Since $\eta<n^{-1}$, for $\varepsilon>0$ sufficiently small we have 
    \begin{align*}
    	\|A-I_n\|_1\le n\eta_\varepsilon<1.
    \end{align*}
    Then $I_n+s(A-I_n)$ is invertible for each $s\in[0,1]$, and $A$ is positive definite. It follows that $x_1,\ldots,x_n$ are linearly independent. Since $\dim\beta=n$, we have 
    \begin{align*}
    	\beta=\Span\{u_1,\ldots,u_n\}.
    \end{align*}
    Thus $hR$ is positive definite on $\beta$.
    \newline 3. Set 
    \begin{align*}
    	(w_1,\ldots,w_n)\ :=(u_1,\ldots,u_n)A^{-\frac{1}{2}}.
    \end{align*}
    Then we have 
    \begin{align*}
    	\left(hR(w_i,w_j)\right)_{i,j=1,\ldots,n}=A^{-\frac{1}{2}}AA^{-\frac{1}{2}}=I_n.
    \end{align*}
    Since there holds
    \begin{align*}
    	\|A^{-\frac{1}{2}}-I_n\|_1=&\|\sum_{k=1}^{+\infty}\left(\begin{array}{c}
    	k	\\-\frac{1}{2}    		
    	\end{array}\right)(A-I_n)^k\|_1\\
        \le&\sum_{k=1}^{+\infty}\left(\begin{array}{c}
        	k	\\-\frac{1}{2}    		
        \end{array}\right)(-\|A-I_n\|_1)^k\|\\
        =&(1-\|A-I_n\|_1)^{-\frac{1}{2}}-1\\
        \le&(1-n\eta_\varepsilon)^{-\frac{1}{2}}-1,
    \end{align*}
    by \eqref{e:norm-u} we have 
    \begin{align*}
    	\|w_i-u_i\|\le& \|A^{-\frac{1}{2}}-I\|_1\max\{\|u_j\|;\;j=1,\ldots,n\}\\
    	\le&((1-n\eta_\varepsilon)^{-\frac{1}{2}}-1)(1+(1+\varepsilon)\delta(\alpha,\beta))\gamma(Q|_\alpha)^{-\frac{1}{2}},\\
    	\|w_i-v_i\|\le& \|w_i-u_i\|+\|u_i-v_i\|\\
    	\le&((1-n\eta_\varepsilon)^{-\frac{1}{2}}(1+(1+\varepsilon)\delta(\alpha,\beta))-1)\gamma(Q|_\alpha)^{-\frac{1}{2}},\\
    	\|w_i\|\le&\|u_i\|+\|w_i-u_i\|\\
    	\le&(1-n\eta_\varepsilon)^{-\frac{1}{2}}(1+(1+\varepsilon)\delta(\alpha,\beta))\gamma(Q|_\alpha)^{-\frac{1}{2}} 
    \end{align*}
    for $i=1,\ldots,n$.
    \newline 4. For each $v\in \mathrm{S}_{\alpha^Q}$, there is a $w\in W$ such that 
    \begin{align*}
    	\|v-w\|\le(1+\varepsilon)\delta(V,W).
    \end{align*}
    Let $i=1,\ldots, n$ be an integer. Note that $Q(v,v_i)=0$. By \eqref{e:c-gap-pair}, Step 1 and Step 3 we have
    \begin{align*}
    	&|R(w,w_i)|=|R(w,w_i)-Q(v,v_i)|\\
    	\le&\delta_c(Q,R)(\|w\|+\|v\|)(\|w_i\|+\|v_i\|)\\
    	&+
    	c\left((\|w\|+\|v\|)\|w_i-v_i\|+\|w-v\|(\|w_i\|+\|v_i\|)\right)\\
    	\le&\gamma(Q|_\alpha)^{-\frac{1}{2}}\delta_c(Q,R)(2+(1+\varepsilon)\delta(V,W))\\
    	&((1-n\eta_\varepsilon)^{-\frac{1}{2}}(1+(1+\varepsilon)\delta(\alpha,\beta))+1)\\
    	&+
    	c\gamma(Q|_\alpha)^{-\frac{1}{2}}(2+(1+\varepsilon)\delta(V,W))\\
    	&((1-n\eta_\varepsilon)^{-\frac{1}{2}}(1+(1+\varepsilon)\delta(\alpha,\beta))-1)\\
    	&+c(1+\varepsilon)\gamma(Q|_\alpha)^{-\frac{1}{2}}\delta(V,W)\\
    	&((1-n\eta_\varepsilon)^{-\frac{1}{2}}(1+(1+\varepsilon)\delta(\alpha,\beta)) +1)\\
    	:\ =&C_\varepsilon.
    \end{align*}
    Define $p_2\in\Bb(W)$ by
    \begin{align*}
    	p_2(x)=x-\sum_{i=1}^{n}hR(x,w_i)w_i,\quad\forall x\in W.
    \end{align*}
    Then we have $p_2(w)\in\beta^R$, and 
    \begin{align*}
    	&\|v-p_2(w)\|\le\|v-w\|+\sum_{i=1}^{n}|R(w,w_i)|\|w_i\|\\
    	\le&(1+\varepsilon)\delta(V,W)
    	+nC_\varepsilon((1-n\eta_\varepsilon)^{-\frac{1}{2}})(1+(1+\varepsilon)\delta(\alpha,\beta))\gamma(Q|_\alpha)^{-\frac{1}{2}} .
    \end{align*} 
    By Definition \ref{d:closed-distance} we have
    \begin{align*}
    	\delta(\alpha^Q,&\beta^R)\le (1+\varepsilon)\delta(V,W)\\
    	&+nC_\varepsilon(1-n\eta_\varepsilon)^{-\frac{1}{2}}(1+(1+\varepsilon)\delta(\alpha,\beta))\gamma(Q|_\alpha)^{-\frac{1}{2}}.
    \end{align*}    
    On letting $\varepsilon\to 0$, we obtain
    \eqref{e:annihilator-close}.  
    \newline 5. For each $\nu\ :=(a_1,\ldots,a_n)\in\KK^n\setminus\{0\}$,
    we  denote by
    \begin{align*}
    	v(\nu)\ :=\sum_{i=1}^{n}a_iv_i,\quad w(\nu)\ :=\sum_{i=1}^{n}a_iw_i.
    \end{align*}
    Then we have
    \begin{align*}
    	hQ(v(\nu),v(\nu))=&hR(w(\nu),w(\nu))\\
    	=&\sum_{i=1}^{n}|a_i|^2\ge\gamma(Q|_\alpha)\|v(\nu)\|^2,\\
    	\|w(\nu)\|\le&\sum_{i=1}^{n}|a_i|\|w_i\|\\
    	\le &\sum_{i=1}^{n}|a_i|(1-n\eta_\varepsilon)^{-\frac{1}{2}}(1+(1+\varepsilon)\delta(\alpha,\beta))\gamma(Q|_\alpha)^{-\frac{1}{2}}\\
    	\le& n^{\frac{1}{2}}(1-n\eta_\varepsilon)^{-\frac{1}{2}}(1+(1+\varepsilon)\delta(\alpha,\beta))\\
    	&\gamma(Q|_\alpha)^{-\frac{1}{2}}(hR(w(\nu),w(\nu)))^\frac{1}{2}.
    \end{align*}
    By definition we have
    \begin{align*}
    	\gamma(R|_\beta)\ge n^{-1}(1-n\eta_\varepsilon)(1+(1+\varepsilon)\delta(\alpha,\beta))^{-1}\gamma(Q|_\alpha).
    \end{align*}
    On letting $\varepsilon\to 0$, we obtain
    \eqref{e:annihilator-close-gamma}.  
\end{proof}

We need the following special case.

\begin{proposition}\label{p:perturbed-Morse-index-definite}
	Let $X$ be a Banach space with two bounded symmetric pairs $(Q,V)$ and $(R,W)$. Let $h\in\{1,-1\}$ and $c\ge 0$ be two real numbers. Assume that  $hQ$ is positive semi-definite and $\gamma(Q)>0$. Let $V_0\subset V^R$ and $W_0\subset W^R$ be two closed linear subspaces of $X$ with finite $n\ :=\dim V^Q/V_0$. Assume that 
	\begin{align}\label{e:perturbed-Morse-index-definite-1}
		\theta\ :=\delta(W,V)+\delta(V_0,W_0)+\delta(W,V)\delta(V_0,W_0)<\frac{1}{2^n(n+1)}.
	\end{align}
	
	Set 
	\begin{align*}
		\delta\ :&=\frac{2^n(n+1)\delta(W,V)}{1-2^n(n+1)\delta(W,V)},\\
		d\ :&=(\gamma(Q))^{-1}, 
		e\ :=2\delta_c(Q,R), \quad 
		f\ :=2\delta_c(Q,R)+c.		
	\end{align*}
	Assume that 
	\begin{align}\label{e:perturbed-Morse-index-definite-2}
		\left(d(2+\delta)(e+f\delta)\right)^{\frac{1}{2}}<\frac{1-2^n(n+1)\theta}{2^n(n+1)(1+\theta)}.		
	\end{align} 
	Then we have 
	\begin{align}\label{e:perturbed-augmented-Morse-index-special}
		m^-(hR)+\dim W^R/W_0\le n.
	\end{align}		
\end{proposition}

\begin{proof} 		
	1. Assume that $m^-(hR)+\dim W^R/W_0>n$. By \cite[Lemma III.1.12]{Ka95}, there is a linearly independent subset $\{u_1+W_0,\ldots,u_{n+1}+W_0\}
	$ of $W/W_0$ such that $(h R(u_i,u_j))_{i,j=1,\ldots,n}\le 0$ and 
	\begin{align*}
		\dist(u_i,\Span\{u_0,\ldots,u_{i-1}\}+W_0)=\dist(u_i,W_0),\quad\text {for }i=1,\ldots,n.
	\end{align*}
	For each $\varepsilon>0$ and $i=1,\ldots,n$, there is a $z_i\in W_0$ such that 
	\begin{align*}
		\dist(u_i,W_0)>(1-\varepsilon)\|u_i-z_i\|>0.
	\end{align*}
	Set 
	\begin{align*}
		w_i\ :&=\frac{u_i-z_i}{\|u_i-z_i\|},\\
		W_i\ :&=\Span\{w_1,\ldots,w_i\}	
	\end{align*} 
	for $i=1,\ldots,n+1$.
	
	Then the subset $\{w_1+W_0,\ldots,w_{n+1}+W_0\}
	$ of $W/W_0$ is linearly independent and we have  
	\begin{align*}
		&(h R(w_i,w_j))_{i,j=1,\ldots,n}\le 0,\quad\|w_k\|=1,\\
		&\dist(w_k,W_{k-1}+W_0)=\dist(w_k,W_0)\\
		&=\dist(\frac{u_k}{\|u_k-z_k\|},W_0)>1-\varepsilon>0
	\end{align*}
	for $k=1,\ldots,n+1$.
	\newline 2. For each $\varepsilon>0$ and $i=1,\ldots,n+1$, there is an $v_i\in V$ such that
	\begin{align*}
		\|v_i-w_i\|<\delta(W,V)+\varepsilon.
	\end{align*}    
	Set 
	\begin{align*}
		V_i\ :&=\Span\{v_1,\ldots,v_i\}		
	\end{align*} 
	for $i=1,\ldots,n+1$. Set
	\begin{align*}
		\Delta_\varepsilon\ :=\frac{(1-\varepsilon)^{n+1}}{(n+1)(2-\varepsilon)^n}.
	\end{align*}
	By Lemma \ref{l:gap-finite-dimension}, for sufficiently small $\varepsilon>0$ we have
	\begin{align*}
		\dim (V_{n+1}+V_0)/V_{n+1}&=\dim (W_{n+1}+W_0)/W_{n+1}=n+1,\\
		\gamma(V_{n+1},V_0)&\ge\frac{\Delta_\varepsilon-(\delta(W,V)+\varepsilon)(1+\delta(V_0,W_0))-\delta(V_0,W_0)}{(1+\delta(W,V)+\varepsilon)(1+\delta(V_0,W_0))},
		\\		
		\hat\delta(V_{n+1},W_{n+1})&\le \frac{\delta(W,V)+\varepsilon}{\Delta_\varepsilon-\delta(W,V)-\varepsilon}.
	\end{align*}
	\newline 3. For each $v\in V_{n+1}$ with $\|v\|=1$, there is a $w\in W_{n+1}$ with
	\begin{align*}
		\|v-w\|\le\delta_\varepsilon\ := \frac{\delta(W,V)+2\varepsilon}{\Delta_\varepsilon-\delta(W,V)-\varepsilon}. 
	\end{align*}
	By \eqref{e:reduced-minimum-modulus-Q} and \eqref{e:c-gap-pair} we have 
	\begin{align*}
		&\dist(v,V^Q)\le(h dQ(v,v))^{\frac{1}{2}}\\
		&\le (d\delta_c(Q,R)(\|v\|+\|w\|)^2+2cd(\|v\|+\|w\|)\|v-w\|)^{\frac{1}{2}}\\
		&\le\left(d(2+\|v-w\|)(e+f\|v-w\|)\right)^{\frac{1}{2}}.
	\end{align*}
	Then we have
	\begin{align*}
		\delta(V_{n+1},V^Q)\le\left(d(2+\delta_\varepsilon)(e+f\delta_\varepsilon)\right)^{\frac{1}{2}}.
	\end{align*}
	\newline 4.
	By \cite[Lemma A.3.1]{BoZh18}, \cite[(1.4.1)]{Neu68}, Step 2 and Step 3, for sufficiently small $\varepsilon>0$ we have
	\begin{align*}
		\delta((V_{n+1}+V_0)/V_0,V^Q/V_0)&= \delta(V_{n+1}+V_0,V^Q)\\
		&\le\delta(V_{n+1},V_Q)\gamma(V_{n+1},V_0)^{-1}<1.
	\end{align*}
	By Step 2 and  \cite[Corollary IV.2.6]{Ka95}, we have
	\begin{align*}
		n+1=\dim (V_{n+1}+V_0)/V_0\le \dim V^Q/V_0=n.
	\end{align*}
	This is a contradiction.
\end{proof}

Now we are able to prove Theorem \ref{t:perturbed-augmented-Morse-index}.

\begin{proof}[Proof of Theorem \ref{t:perturbed-augmented-Morse-index}] 
	If $\dim\alpha=+\infty$ we are done. Assume that $\dim\alpha<+\infty$. Then there is a $Q$-orthogonal direct sum decomposition
	\begin{align*}
		\alpha=\alpha_1\oplus\alpha_0
	\end{align*}
	with $\alpha_0=\alpha^Q\cap\alpha$. 
	Set $k\ :=\dim\alpha_1$ and $l\ :=\dim\alpha_0$. By Lemma \ref{l:existence-finite-dimension-embeding-close}, 
	for each $\varepsilon\in(0,\frac{1}{2^{k-1}k}-\delta(V,W))$, there is a linear subspace $\beta_1$ of $W$ such that there hold $\beta_1=k$ and
	\begin{align}\label{e:morse-existence-finite-dimension-embeding-close}
		\hat\delta(\alpha_1,\beta_1)\le\frac{2^{k-1}k(\delta(V,W)+\varepsilon)}{1-2^{k-1}k(\delta(V,W)+\varepsilon)}.
	\end{align}
    Since $\delta(V,W)$ is sufficiently small, $\hat\delta(\alpha_1,\beta_1)$ is also sufficiently small. By \cite[Corollary IV.2.6]{Ka95} we have $\dim\alpha_1=\dim\beta_1$. By \cite[Lemma 2.15.e]{LhcZ24} we have $\delta_c(Q|{\alpha_1},R|_{\beta_1})\le\delta_c(Q,R)$. By \cite[Proposition 2.21]{LhcZ24} we have 
    \begin{align*}
    	\dim\beta_1\ge m^+(hR|_{\beta_1},\beta_1)\ge\dim\alpha_1=\dim\beta_1.
    \end{align*}
    Then $hR$ is positive definite on $\beta_1$.
	By Lemma \ref{l:annihilator-close}, $\hat\delta(\alpha_1^Q,\beta_1^R)$ is sufficiently small. Since $W^R\subset\beta_1^R$, by \cite[Lemma 1.1.2.c, Lemma A.1.1]{BoZh18} we have
	\begin{align*}
		\beta_1^{RR}\cap\beta_1^R=(\beta_1+\beta_1^R)^R=(\beta_1+W^R)\cap\beta_1^R=\beta_1\cap\beta_1^R+W^R=W^R.
	\end{align*}
	By \eqref{e:V-a-b} we have $\alpha_1^Q=\alpha_0\oplus\beta$. By \cite[Lemma A.1.1]{BoZh18} we have 
	\begin{align*}
		(\alpha_1^Q)^{Q|\alpha_1}=&(\alpha_0\oplus\beta)^Q\cap(\alpha_0\oplus\beta)=(\alpha_0^Q\cap\beta^Q)\cap(\alpha_0\oplus\beta)\\
		=&\beta^Q\cap(\alpha_0\oplus\beta)
		=\alpha_0\oplus\beta^Q\cap\beta\\
		=&\alpha_0\oplus V_0.
	\end{align*}
	By Proposition \ref{p:perturbed-Morse-index-definite}, we have 
	\begin{align*}
		m^{+}(hR)+\dim W^R/W_0=&k+m^{+}(hR|_{\beta_1^R})+\dim W^R/W_0\\
		\le &k+\dim\alpha_0=\dim\alpha.
	\end{align*}
	
\end{proof}


\bibliography{Hamiltonian}
\bibliographystyle{plain}

\end{document}